\documentclass[letterpaper,12pt,peerreviewca,draftcls]{IEEEtran}
\usepackage{csmMBFeb07,graphicx,url}
\usepackage{color}
\usepackage{cite}
\usepackage{amsmath,amssymb,theorem}
\usepackage{graphicx}
\usepackage{algorithm,algorithmic}
\usepackage{psfrag}
\usepackage{url}
\usepackage{lipsum}
\usepackage{pbox}
\usepackage{mdframed}
\usepackage{tikz}
\usepackage[utf8]{inputenc}
\usepackage{csquotes}
\usetikzlibrary{automata,positioning}

\newtheorem{theorem}{Theorem}

\newtheorem{remark}[theorem]{Remark}
	
\newtheorem{corollary}[theorem]{Corollary}

\newtheorem{definition}[theorem]{Definition}
\newtheorem{problem}[theorem]{Problem}

\newmdtheoremenv{sidebar}{Sidebar}

\newcommand{\probability}[2]{\operatorname{Pr}(#1| \, #2)}

\newcommand{\real}{{\mathbb{R}}}
\newcommand{\realpositive}{\mathbb{R}_{>0}}
\newcommand{\realnonnegative}{\mathbb{R}_{\ge 0}}

\newcommand{\GG}{{\mathcal{G}}}

\newcommand{\EE}{{{E}}}

\newcommand{\NN}{{\mathcal{N}}}

\newcommand{\VV}{{{V}}}

\newcommand{\diag}[1]{\operatorname{diag}\left( #1\right)}

\renewcommand{\tilde}{\widetilde}

\newcommand{\timestep}{\Delta t}

\renewcommand{\epsilon}{\varepsilon}

\newcommand{\oprocendsymbol}{\hbox{$\bullet$}}
\newcommand{\oprocend}{\relax\ifmmode\else\unskip\hfill\fi\oprocendsymbol}

\newcommand{\until}[1]{\{1,\dots, #1\}}

\newcommand{\setdef}[2]{\{#1 \; | \; #2\}}

\newcommand{\Nin}{\NN^\text{in}}
\newcommand{\Nout}{\NN^\text{out}}

\newcommand{\betaeff}{\beta^\text{eff}}

\title{Analysis and Control of Epidemics \\
{\Large A survey of spreading processes on complex networks}
}
\author{Cameron Nowzari, Victor M. Preciado, and George J. Pappas \\ \today}

\begin{document}

\maketitle
\CSMsetup

\begin{center}
{\color{blue}NOTE: This arXiv version contains a table of contents at the end for convenience.}
\end{center}

This article reviews and presents various solved and open problems in the development, analysis, and control of epidemic models. Proper modeling and analysis of spreading processes has been a longstanding area of research among many different fields including mathematical biology,
physics, computer science, engineering, economics, and the social sciences. 
One of the earliest epidemic models conceived was by Daniel Bernoulli in 1760,
which was motivated by studying the spreading of smallpox~\cite{DB:1760}. In addition to
Bernoulli, there were many different researchers also working on mathematical epidemic
models around this time~\cite{KD-JAPH:02}. These initial models were quite simplistic
and the further development and study of such models dates back 
to the 1900s~\cite{RR-HPH:16,RR-HPH:17,RR-HPH:17b,WOK-AGM:27}, 
where still simple models were studied to provide insight as to how various diseases
can spread through a population. In recent years, we have seen a resurgence of interest
in these problems as the concept of `networks' becomes increasingly prevalent in modeling
many different aspects of our world today. A more comprehensive review of the history of
mathematical epidemiology can be found in~\cite{HWH:00,NTB:75}. 

Despite the study of epidemic models having spanned such a long period of time, it is only
recently that control engineers have entered the scene. Consequently, there is already a vast
body of work dedicated to the development and analysis of epidemic models, but far less
that provide proper insight and machinery on how to effectively \emph{control} these processes. 
The focus of this article is to provide an introductory tutorial on the latter. 
We are interested in presenting
a relatively concise report for new engineers looking to enter the field of spreading
processes on complex networks. This article
presents some of the more well-known and recent results in the literature while also identifying
numerous open problems that can benefit from the
collective knowledge of optimization and control theorists. 

Although this article focuses on the context of epidemics, the same models and tools
we present are directly applicable to a myriad of different
spreading processes on complex networks. Examples include the adoption of an idea 
or rumor through a social network like Twitter,
the consumption of a new product in a marketplace, the risk of receiving a computer virus
through the World Wide Web, and of course the spreading of a disease through a 
population~\cite{SB-VL-YM-MC-DUH:06,DE-JK:10,MEJN:10}.
For this reason, we will often use the terms individuals, people, nodes, and agents
interchangeably. 

We begin the article by introducing and analyzing some classical epidemic models
and their extensions to network models. We then discuss various methods of controlling
these epidemic models and several extensions. After describing the main
shortcomings in the current literature and highlighting some recent preliminary works that are aimed
at improving the current state of the art, we close by providing some intuition into
the current research challenges that need to be addressed in order to fully harness the power
of these works and make a real societal impact.

\section{Modeling and Analysis of Epidemics}


Before jumping into the class of models we study in this article, we must start by
emphasizing that there are an uncountable number of ways to model spreading processes.
The underlying common factor that ties almost all epidemic models together is the existence
of `compartments' in which individuals in a population are divided. The two most common compartments that
exist in essentially every single epidemic model are called `Susceptible' and `Infected'~\cite{WOK-AGM:27,HWH:00,RMA-RMM-BA:92}. In models that contain only these two compartments, a given population is initially divided into
these two compartments. The `Susceptible' compartment ($S$) represents individuals who are healthy, but susceptible to becoming infected. The 'Infected' compartment ($I$) captures individuals who are infected, but
able to recover. From here there are an insurmountable number of ways that the interactions within
the population can be modeled. 

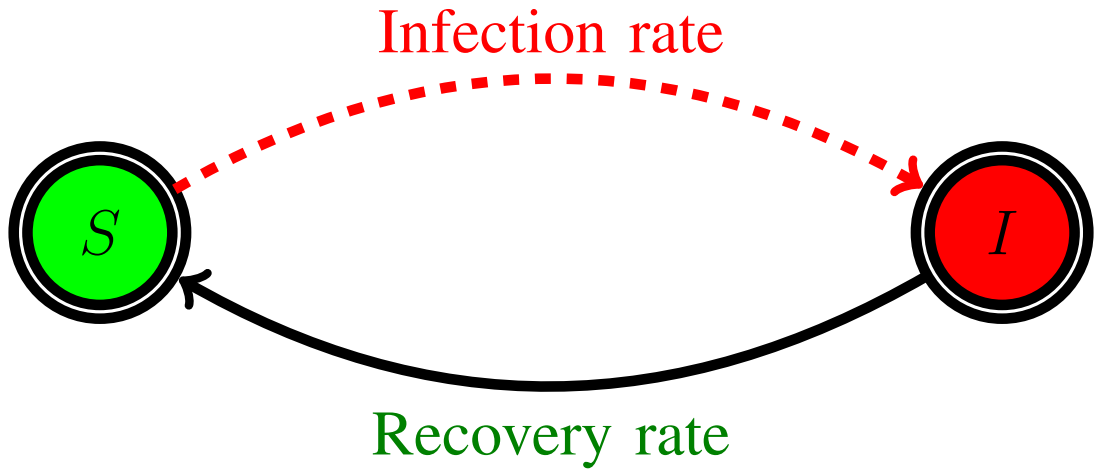
\begin{figure}[htb]
  \centering
  \scalebox{2}{\begin{tikzpicture}[shorten >= 1pt,node distance=2.25cm,auto]
  \pgfsetlinewidth{2pt}
    \node[state, fill=green, accepting] (S)              {$S$};
    \node (G) [right=of S] {};
    \node[state, fill=red, accepting] (I) [right=of G] {$I$};
    
    \path[->] (I) edge [bend left] node {{\color{green!50!black}Recovery rate}} (S);
    
    \path[->,red,dashed] (S) edge [bend left] node{Infection rate} (I);
    \end{tikzpicture}}
    \caption{Two-state Susceptible-Infected-Susceptible (SIS) model. An individual
    in the infected state~$I$ transitions to the healthy or susceptible state~$S$ with
    some recovery rate and from the susceptible state to the infected state with
    some infection rate.}\label{fig:SIS_compartmental}
\end{figure}

Throughout this article, we focus on models where individuals can move from one compartment to another
randomly with some defined rates. For instance, in this two compartment model healthy individuals can randomly transition from $S$ to $I$ with some infection rate that is a result of interactions with infected individuals. Similarly, infected individuals can randomly transition from $I$ to $S$ with some recovery rate that is a result of recovering from the infection. More details on how these rates are defined are provided later.
Figure~\ref{fig:SIS_compartmental} shows the simple interaction described above.

In addition to models with only two compartments, there are many other epidemic models aimed at
capturing various important features of realistic diseases and spreading processes. 
This is often done by adding more compartments,
such as a `Removed' ($R$) compartment representing individuals who are no longer susceptible to the infection. This might refer to a deceased, vaccinated, or immune individual. Other compartments have also been proposed in the literature to study the effect of, for example, an incubation period, partial immunity, or quarantine in the spreading dynamics~\cite{PD-JW:02,SF-EG-VAAJ:10,IZK-JC-MR-PLS:10,NP-DB-BG-AV:11,PP-BC-MA-AP-SM:09,DB-OD-WFG-AP-RV:12,MMH:14}. 

For brevity, we will now focus our attention on two of the oldest epidemic models known
as the Susceptible-Infected-Removed (SIR) and Susceptible-Infected-Susceptible (SIS) 
models~\cite{WOK-AGM:27}. However, we note that
the following exposition can be applied to more general compartmental models with some appropriate changes.

Let $N$ be the total number of individuals in a population. We denote by $X_i(t) \in \{ S, I, R \}$
the state of node $i \in \until{N}$ at time $t$. We collect the state of the entire population
in a state vector $X(t) = (X_1(t), \dots, X_N(t))^T$. The evolution of the states is then described
by a Markov process as follows. An individual $i$ infected at time $t$ recovers to the removed state at a fixed rate $\delta_i > 0$. In other words, if node $i$ is infected at time $t$, the probability that this node 
loses its infection in the time slot $(t,t+\Delta t]$ for small $\Delta t$ is given by
\begin{align*}
\probability{ X_i(t+\Delta t) = R}{ X_i(t) = I } &= \delta_i \Delta t +o(\Delta t), \quad \text{(SIR)} \\
\probability{ X_i(t+\Delta t) = S}{ X_i(t) = I } &= \delta_i \Delta t +o(\Delta t). \quad \hspace*{.4ex} \text{(SIS)}
\end{align*}
The above represents an
\emph{endogenous} transition which occurs internally within each node, independent of the states of other nodes~\cite{BAP-DC-NCV-MF-CF:12}.

\begin{figure}[htb]
  \centering
    \scalebox{2}{{\begin{tikzpicture}[shorten >= 1pt,node distance=2.25cm,auto]
  \pgfsetlinewidth{2pt}
    \node[state, fill=green, accepting] (S)              {$S$};
    \node[state, fill=red, accepting] (I) [right=of S] {$I$};
    \node[state, fill=blue!80!white, accepting] (R) [right=of I] {$R$};
    
    \path[->] (I) edge [bend left] node {{\color{green!50!black}$\delta_i$}} (R);
    
    \path[->,red,dashed] (S) edge [bend left] node{$\beta^\text{eff}_i$} (I);
    \end{tikzpicture}}}
    \caption{Three-state Susceptible-Infected-Removed (SIR) model. An individual~$i$
    in the susceptible state~$S$ can transition to the infected state~$I$ with
    some infection rate~$\beta^\text{eff}_i$ and from the infected state~$I$ to
    the removed state~$R$ with some recovery rate~$\delta_i$.}\label{fig:SIR_compartmental}
\end{figure}
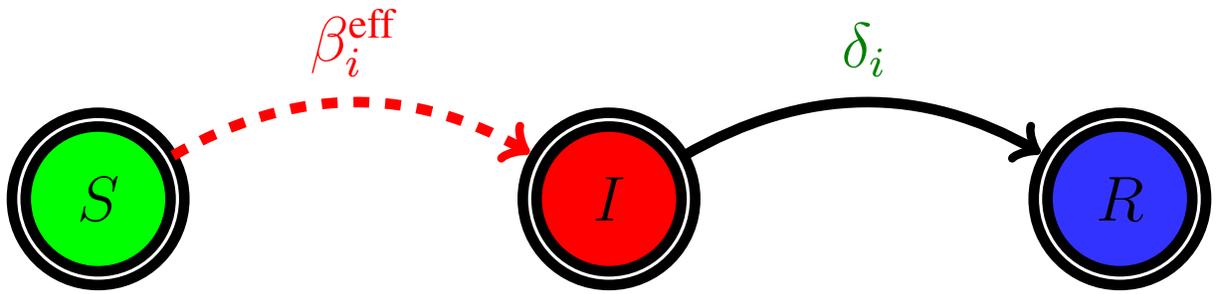

Similarly, an individual~$i$ that is susceptible at time~$t$ becomes infected
at a rate $\beta^\text{eff}_i$ that depends
on the state of the entire population $X(t)$. This is known
as an \emph{exogenous} transition because it is influenced by factors external
to the node itself. We discuss this at length in the sections to come. 
Figure~\ref{fig:SIR_compartmental} shows the simple interaction described above.

\begin{remark}[Other spreading models]
{\rm
We note here that this article excludes chain binomial models (e.g., Reed-Frost model) and
other similar types of models from percolation theory. Depending on the application at hand,
the model for the spreading dynamics can vary. The main difference in the models we consider and models like Reed-Frost are that we want to allow infected individuals to continuously try to infect healthy ones. In the Reed-Frost
model, an infected person only has one chance of infecting a healthy person. However, when thinking of
a virus like the flu, a healthy person is continuously in danger of becoming sick when in contact with
an infected individual, rather than a one-time chance. For instance, the Reed-Frost model might be
more suitable for modeling the spreading of an email virus rather than an infectious biological
disease. The interested reader should see~\cite{NTB:75,TB:10} for further details.

We would also like to briefly point out that this other community is indeed active and working on similar
types of problems as the ones we will highlight in this article. Many works exist on forecasting
the cascading effects of a single infection or failure on a network~\cite{MG-WG-DT:03,LMS-CPW-IMS-JK:02}
and how they can be mitigated through vaccination~\cite{EK-JCM:11}. 
Conversely, one may be more interested in finding the most influential nodes or where to start
an infection in a network to reach as many people as possible~\cite{DK-JK-ET:03,JK:07}. 
This is often referred to as a seeding problem. Further extensions study attack and vaccination strategies on these models~\cite{ML:09} and even cases in which there are multiple contagions on multiple networks~\cite{OY-VG:12}. \oprocend
}
\end{remark}

\subsection{Classical models}

\begin{figure}[htb]
  \centering
  {\includegraphics[width=.85\linewidth]{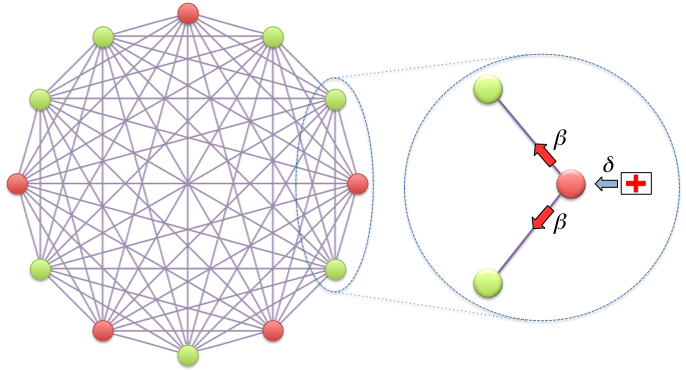}}
  \caption{Population dynamics of the two-state SIS model. These models assume a well-mixed population,
  meaning that each individual in the population is equally likely to contract a disease from anyone
  else in the population. An infected individual (red) naturally recovers at a rate~$\delta > 0$,
  depicted by the red cross.
  A healthy individual (blue) if affected by each infected
  individual in the population with rate~$\beta$, depicted by the red arrows.
  }
  \label{fig:SIS_population}
\end{figure}

Based on the above discussion, the dynamics of the SIR model
is described by a $3^N$-dimensional Markov process. The exponential size of the state space makes this model very hard to analyze. One standard method to simplify the analysis is to consider the evolution of the total number
of healthy and infected individuals rather than the state of each individual separately.
This is commonly referred to as \emph{population dynamics}~\cite{RMA:82,FB-CCC:12}.
Furthermore, the recovery and infection rates are often assumed to be the same for all individuals; that is, $\delta_i=\delta$ and $\beta^\text{eff}_i=\beta^\text{eff}$ for all $i$. The standard population dynamics assumes a 
\emph{well-mixed population} which means all individuals affect and are affected
by all other individuals equally. Figure~\ref{fig:SIS_population} shows the
described interactions of this well-mixed population.

\subsubsection{Stochastic population models}

The SIR population model is described as follows. 
Letting $N^I(t), N^R(t) \in \{0, 1, \dots, N\}$ be the number of infected and removed individuals 
at some time $t$, respectively, the number of susceptible individuals is necessarily given
by~$N^S(t) = 1 - N^I(t) - N^R(t)$. A common choice for the 
infection rate is given by $\betaeff = \beta N^I N^S$~\cite{NTB:57,KD:67,HWH:00}
for some $\beta > 0$, known as the mass action law. In other words, the rate at which the total number of susceptible
individuals become infected is proportional to the product of the number of susceptible
and infected individuals in the population. The state at some time $t + \timestep$
is then given by
\begin{align}\label{eq:population_stochastic_SIR}
(N^I, N^R) \rightarrow \begin{cases} (N^I + 1, N^R) & \text{ with probability } \beta N^I N^S \timestep + o(\timestep), \\
(N^I - 1, N^R - 1) & \text{ with probability } \delta N^I \timestep + o(\timestep), \\
(N^S, N^R) & \text{ with probability } 1 - (\beta N^I N^S + \delta N^I) \timestep + o(\timestep).
\end{cases}
\end{align}
For the SIS model we simply force $N^R = 0$ at all times which simplifies this to
\begin{align}\label{eq:population_stochastic_SIS}
N^I \rightarrow \begin{cases} N^I+1 & \text{ with probability } \beta N^I N^S \timestep + o(\timestep), \\
N^I-1 & \text{ with probability } \delta N^I \timestep + o(\timestep), \\
N^I & \text{ with probability } 1 - (\beta N^I N^S + \delta N^I) \timestep + o(\timestep).
\end{cases}
\end{align}

Removing the explicit definition of time, the SIS process
can then be seen as a random walk on a line~\cite{KD:67,AA:73,AA:74,RJK-CL:89} (a similar Markov chain can be described for the SIR model)

\begin{align}\label{eq:population_random_walk}
N^I \rightarrow N^I + 1 \text{ with probability } \frac{\beta(N-N^I)}{\beta(N-N^I) + \delta} ,\\
N^I \rightarrow N^I - 1 \text{ with probability } \frac{\delta}{\beta(N-N^I) + \delta} . \notag
\end{align}

An important observation about this model~\eqref{eq:population_random_walk} 
is that it is a Markov chain with a single
absorbing state~$N^I = 0$ in which all agents are healthy. In other words, once the entire population is healthy, the infection cannot suddenly reemerge. It is known from the theory of Markov chains that given enough time, the infection will eventually die out with probability~1
(see~\cite{PB:99} for a review of Markov chains and relevant properties).
Thus, the study of these systems is often interested in answering the question of \emph{when} or \emph{how quickly} the infection
will die out. We comment on this later, in Remark~\ref{re:population_comparison}.

To further simplify the problem, many works consider a deterministic approximation of these
stochastic dynamics. In fact, the simpler deterministic dynamics we introduce next predate
the introduction of the stochastic model above~\cite{WOK-AGM:27}. 

\subsubsection{Deterministic population models}

The models we present next are perhaps the two most studied epidemic models in the literature 
and is covered in a large number of different books~\cite{WOK-AGM:27,NTB:75,RMA-RMM-BA:92,DJD-JMG:99,OD-JAPH:00,HA-TB:00,MJK-PR:07,PM-SR:08,PVM:09,MD-LM:10,DE-JK:10,MEJN:10,OD-HH-TB:12}. These books also
discuss a large variety of extensions
including more complicated disease models that have more than two states, modeling for birth and mortality
rates, different types of infection rates, and different categories for each disease state;
for example, based on age or sex. We only present the most basic models here.


Assuming a large population size $N$, we define $p^I = \frac{N^I}{N}$ and $p^S = \frac{N-N^I-N^R}{N}$ 
as the fractions of infected and susceptible individuals, respectively.
Then, we can write the deterministic SIR version of~\eqref{eq:population_stochastic_SIR} as
\begin{align}\label{eq:population_deterministic_whole_SIR}
\dot{p}^S &= -\beta p^I p^S, \\
\dot{p}^I &= \beta p^I p^S - \delta p^I , \notag
\end{align}
and the deterministic SIS version of~\eqref{eq:population_stochastic_SIS} as
\begin{align}\label{eq:population_deterministic_whole_SIS}
\dot{p}^S &= -\beta p^I p^S + \delta p^I, \\
\dot{p}^I &= \beta p^I p^S - \delta p^I . \notag
\end{align}
These are derived by leveraging Kurtz' theorem while assuming~$N$ to be very large~\cite{MD-LM:10}.

For simplicity, we will continue the analysis only for the SIS model but note that similar analysis can be
done for the SIR model as well. Since the population size $N$ is fixed and $N^R = 0$, we have that $p^S = 1 - p^I$, and the above equations~\eqref{eq:population_deterministic_whole_SIS} are redundant. Hence, they can be simplified to
\begin{align}\label{eq:population_deterministic}
\dot{p}^I &= \beta p^I (1 - p^I) - \delta p^I . 
\end{align}
Given an initial condition $p^I(0)$, this equation can be analytically solved~\cite{GHW-MD:71,NTB:75,HWH:76}. The solution is given by
\begin{align*}
p^I(t) = \left\lbrace \begin{array}{lc} \frac{e^{(\beta-\delta)t}}{ \frac{ \beta(e^{(\beta-\delta)t}-1)}{\beta-\delta} + \frac{1}{p^I(0)}} , & \beta \neq \delta , \\
\frac{1}{ \beta t + \frac{1}{p^I(0)}} , & \beta = \delta . \end{array} \right.
\end{align*}
Given the exact solution of $p^I(t)$, we are able to characterize its equilibrium points in the following result.

\begin{theorem}[Solutions to deterministic population model]\label{th:population_threshold}
The solution of $p^I(t)$ approaches $1 - \frac{\delta}{\beta}$ as $t \rightarrow \infty$ for $\beta > \delta$,
and $0$ as $t \rightarrow \infty$ for $\beta \leq \delta$.
\end{theorem}

\begin{remark}[Deterministic vs stochastic population models]\label{re:population_comparison}
{\rm
It is important to note here that the deterministic models are only approximations of the stochastic models.
A natural question to ask then is the following: What can the threshold result of Theorem~\ref{th:population_threshold} tell us about the stochastic model in~\eqref{eq:population_random_walk}? The first thing to note is that in the stochastic model, given
enough time, the system will reach the disease-free state with probability 1. However, from Theorem~\ref{th:population_threshold} we see that
for $\beta > \delta$, the deterministic model will converge to an endemic equilibrium, meaning the disease
never dies out. Thus, rather than studying the equilibrium values of the two models, the authors in~\cite{GHW-MD:71,RHN:82,RJK-CL:89} look at the expected time~$E[T]$ for the stochastic model to reach the disease-free equilibrium. Interestingly,
they are able to show that for $\beta < \delta$, the expected time $E[T]$ is upper-bounded by $\frac{N \beta}{\delta}$. On the other hand, when $\beta > \delta$, the expected time $E[T]$ grows exponentially with $N$. \oprocend
}
\end{remark}

The analysis of the deterministic model results in a very precise threshold result
that translates directly to the stochastic model as discussed in Remark~\ref{re:population_comparison}. 
Threshold conditions are often given in terms of a \emph{reproduction number} $R_0$, which is the expected number of individuals a single infected individual will infect~\cite{OD-JAPH-JAJM:90,HWH:00} over the course of its infection period. In other words, given a fully healthy population, if we infect person~$i$ at random, $R_0$ is the expected number of other individuals that will become infected over the course of agent~$i$'s infection. This is a very
useful metric with a critical value of $R_0 = 1$. When $R_0 < 1$ the disease does not spread quickly
enough, resulting in a decay in the number of infected individuals (in expectation). On the other
hand, when $R_0 > 1$ the infected population grows over time (in expectation)~\cite{DE-JK:10}. 
In the simple model considered above, the reproduction number is given by $\frac{\beta}{\delta}$. 
Furthermore, the exact solutions and asymptotic behavior of the system are easy to obtain. 

The reproduction number is an important parameter that epidemiologists are interested in identifying
for different diseases and environments~\cite{KD:93} as it is a single number that can predict whether
a certain outbreak of a disease will become an epidemic or die out on its own.
Of course, the problem is that computing~$R_0$ in general is not trivial, as there is no database for
things like infection rates and recovery rates for different diseases.

The main drawback of these population models is that they are very crude model
derived by making many simplifying assumptions including (\emph{i}) a homogeneous incidence rate $\betaeff$ and
recovery rate $\delta$ for all individuals, (\emph{ii}) a two-state model, (\emph{iii}) a constant population size, 
and (\emph{iv}) a well-mixed population (or a contact network that is a complete graph). 
These drawbacks were very evident when scientists attempted to estimate the reproduction number 
of SARS in China in 2002-2003 and grossly overestimated it. This led to SARS scares making global headlines
which eventually fizzled out because the actual reproduction number was far less than estimated due
to the crude population models. More details on how this occurred can be found in~\cite{LAM:07}, but
the upshot is that more refined models are needed. 

%

%
%
%
%
%

\subsection{Network models}

In order to create more refined epidemic models, it is clear that we cannot simply lump an entire population
into two compartments defined by a single number. Ideally, we would be able to model the 
states of all~$N$ individuals independently and allow for arbitrary interactions among them. 
Not surprisingly, this is not a trivial task. 

In this section we are interested in spreading processes on a given, arbitrarily topology.
Before jumping into the models we are interested in, we should mention that we are skipping over a plethora
of work that has been dedicated to extending the population models to \emph{structured} network models. 
More specifically, before jumping to completely arbitrary networks, there is a large body
of work that studies various, specific structures. For instance, some works study how a disease
spreads on a two-dimensional lattice or star graph~\cite{MEJN:02,KO-EGN-AK:13,AMM-AS:13}. Others consider more complex interconnection patterns, such as power-law and small-world networks, that still have some exploitable structure. In this context, a common method to analyze these networks is to assume that nodes are infected at a rate proportional to the number of neighbors they have~\cite{MB-AB-RPS-AV:05,RPS-AV:01,REK-PS-CS-MY:09,YM-RPS-AV:02,MY-RK-CS:11,MY-CS:11}. These methods are justified depending on
the assumptions enforced on the network topology. A review of these types of models 
can be found in~\cite{RPS-CC-PVM-AV:14}. In what follows, we present 
epidemiological models on arbitrary networks.

\subsubsection{Stochastic network models}

Here we describe an epidemic model described as a continuous-time networked Markov model. Consider a network of $N$ nodes represented by a connected, undirected graph $\GG = (V,E)$ where $V$ is the set of nodes and $E \subset V \times V$ is the set of edges. We define the matrix $A \in \realnonnegative^{N \times N}$ as the adjacency matrix of the graph, defined component-wise as
$a_{ij} = 1$ if node $i$ can be directly affected by node~$j$, and $a_{ij} = 0$ otherwise.
See ``Graph Theory'' for further details.

\begin{figure}[htb]
  \centering
  {\includegraphics[width=.85\linewidth]{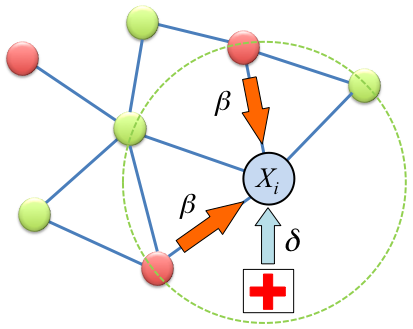}}
  \caption{Network dynamics of the two-state SIS model. A node~$i$ has a natural
  recovery rate $\delta$, depicted by the red cross, 
  at which it transitions from the infected state~$I$
  to the susceptible state~$S$ and is affected by each infected neighbor~$j$
  with rate~$\beta$, depicted by the red arrows.}  \label{fig:SIS_network}
\end{figure}

We let $X_i(t)$ denote the state of node $i$ at time $t$, where $X_i(t) = 1$ indicates that $i$ is infected and $X_i(t) = 0$ indicates that $i$ is healthy at time $t$. Infected nodes can transmit the disease
to its neighbors in the graph~$\GG$ with rate $\beta > 0$. Simultaneously, infected nodes recover from the disease with rate $\delta > 0$. Figure~\ref{fig:SIS_network} shows the described interactions on an arbitrary network. One can model this spreading process using the Markov process
\begin{align}\label{eq:stochastic_network_sis}
\begin{array}{ll} X_i : 0 \rightarrow 1 & \text{with rate } \beta \sum_{j \in \NN_i} X_j ,\\
                  X_i : 1 \rightarrow 0 & \text{with rate } \delta . \end{array} 
\end{align}
Notice that there exists one absorbing
state in this Markov process (corresponding to the disease-free equilibrium) that can be reached from any state
$X(t) = [X_1(t), \dots, X_N(t)]^T$. This implies that, regardless of the initial condition $X(0)$, the epidemic will eventually die out in finite time with probability 1.

An interesting measure of the virality of a spreading process is the expected time $E[T]$ it takes for the epidemic to die out. In~\cite{AG-LM-DT:05} and~\cite{YW-DC-CW-CF:03}, we find the following threshold conditions in terms of the infection strength, defined as $\tau = \frac{\beta}{\delta}$.

\begin{theorem}[Threshold for sublinear expected time to extinction]\label{th:SIS-threshold-homogeneous}
If
\begin{align*}
\tau < \frac{1}{\lambda_\text{max}(A)},
\end{align*}
where $\lambda_\text{max}(A)$ is the maximum real eigenvalue of~$A$, then
\begin{align*}
E[T] \leq \frac{\log N + 1}{\delta - \beta \lambda_\text{max}(A)}
\end{align*}
for any initial condition $X(0)$. 
\end{theorem}

Note that Theorem~\ref{th:SIS-threshold-homogeneous} only provides a sufficient condition
for `fast' extinction of a disease. There have been many efforts to determine
whether this condition is also necessary, but at the time of writing this remains
an open question on general graphs. The works~\cite{MD-LM:10,TM-JCM-DV-QY:12,PVM:13} show that
there exists some critical value $\tau_c$ of the infection strength for which
the expected time to extinction grows exponentially with $N$ when $\tau \geq \tau_c$.
The following result formalizes this statement and provides a lower bound 
on the critical values~\cite{PVM:14,PVM-FDH-CS:14}; however, we should also note that
stronger statements exist when considering graphs with a fixed 
structure (e.g., lattice, star)~\cite{MD-LM:10}.

\begin{theorem}[Threshold for exponential expected time to extinction]\label{th:SIS-threshold-homogeneous-upper}
There exists 
\begin{align*}
\tau_c \geq \frac{1}{\lambda_\text{max}(A)}
\end{align*}
such that for $\tau > \tau_c$, the expected time to extinction $E[T] = O\left( e^{kN} \right)$,
where $k$ depends on $\tau$ and the structure of the graph~$\GG$.
\end{theorem}

The maximum eigenvalue $\lambda_\text{max}(A)$ of an adjacency matrix is a parameter that
captures how `tightly connected' the graph is. In general, more connections means 
a large~$\lambda_\text{max}(A)$. Intuitively, the results of Theorems~\ref{th:SIS-threshold-homogeneous}
and~\ref{th:SIS-threshold-homogeneous-upper} are saying that the more tightly connected the
graph is, the easier it is for a disease to spread. 

It is worth mentioning that although the result of Theorem~\ref{th:SIS-threshold-homogeneous}
provides an upper bound on the \emph{expectation} of the extinction time, the possibility
of a persisting epidemic is not ruled out. For example, it has been shown for star graphs
that regardless of the infection strength $\tau$, there is a positive probability
that the time to extinction is super-polynomial in the number 
of nodes~\cite{CB-JC-AG-AS:10,NB-CB-JTC-AS:05,AG-LM-DT:05}.
Furthermore, for high-degree or scale-free networks (such as preferential attachment~\cite{NB-CB-JTC-AS:05}
or power-law configuration model graphs~\cite{AG-LM-DT:05}), it has been shown that this threshold
goes to zero as the number of nodes increases~\cite{RPS-AV:01b} because the maximum eigenvalue
grows unbounded with $N$. 


\subsubsection{Deterministic network models}

Here we describe the deterministic model of the SIS dynamics over arbitrary 
networks~\cite{DC-YW-CW-JL-CF:08,PVM-JO-RK:09,BAP-DC-NCV-MF-CF:12,PVM:11,CL-RB-PVM:12,HJA-BH:13}.
We begin by assuming homogeneous recovery and infection rates, although this assumption will be relaxed in the
following section.
The natural recovery rate of each node is given by $\delta > 0$ and
the infection rate at which a node is affected by infected neighboring nodes is $\beta > 0$.
The dynamics of the spread is described by the set of ordinary differential equations
\begin{align}\label{eq:network_sis}
\dot{p}_i &= - \delta p_i + \sum_{j=1}^N a_{ij} \beta p_j (1 - p_i),
\end{align}
where $p_i(t) \in [0,1]$ describes the (approximated) probability that 
an individual $i$ is infected at time $t$. 
See ``Networked Mean-Field Approximations'' for further details.
This variable has another interesting interpretation in the context of metapopulation models. In a metapopulation model, each node does not represent an individual, but a large subpopulation (such as an entire district or city). In this context, $p_i$ can be interpreted as the fraction of the $i$-th subpopulation that is infected. See ``Meta-Population Models'' for further details.

As with all other epidemic models, we see that the disease-free equilibrium $p_i = 0$ for all $i \in \until{N}$ 
is a trivial equilibrium of the dynamics. We are now interested in finding conditions such that this equilibrium
is globally asymptotically stable. Letting $p = (p_1, \dots, p_N)^T$ 
and recalling the infection strength $\tau = \frac{\beta}{\delta}$, 
the following result from~\cite{AL-JAY:76,HJA-BH:13,AK-TB-BG:14,AK-TB-BG:14b} characterizes
the convergence properties of these dynamics.

\begin{theorem}[Threshold condition for networks]\label{th:network_homo_threshold}
Given the dynamics~\eqref{eq:network_sis} for any $p(0) \neq 0$, the equilibrium $p^* = 0$ is
globally asymptotically stable if and only if $\tau \leq \frac{1}{\lambda_\text{max}(A)}$. Furthermore,
for $\tau > \frac{1}{\lambda_\text{max}(A)}$, there exists $p^{**} \in \real_{(0,1)}^N$ such that
$p^{**}$ is globally asymptotically stable. 
\end{theorem}

\begin{remark}[Deterministic vs stochastic network models]\label{re:network}
{\rm
Similar to our discussion in Remark~\ref{re:population_comparison},
there is a connection between the deterministic result in Theorem~\ref{th:network_homo_threshold}
and the stochastic result in Theorem~\ref{th:SIS-threshold-homogeneous}.
Since $X = 0$ is an absorbing state, the stochastic dynamics will eventually reach the disease-free state with probability 1. However, Theorem~\ref{th:network_homo_threshold} claims that
for $\beta \lambda_\text{max}(A) > \delta$ the deterministic model will converge to an endemic equilibrium, meaning the disease never dies out. To resolve this apparent contradiction, we again focus our attention on the expected time~$E[T]$ for the stochastic model to reach the disease-free equilibrium.
Remarkably, Theorem~\ref{th:SIS-threshold-homogeneous} provides a sufficient condition
for a disease to quickly die out that is in agreement with the threshold result of
Theorem~\ref{th:network_homo_threshold}. However, as suggested by Theorem~\ref{th:SIS-threshold-homogeneous-upper}, it has not yet been shown whether the same threshold condition holds for persistence of the disease in the stochastic model. \oprocend
}
\end{remark}

A major drawback of the above model is that it assumes a constant infection rate~$\beta$ and recovery rate~$\delta$ for all individuals. To further refine the model, we are interested in allowing different recovery rates for each person and different infection rates for each type of contact. This allows for a much more general model that can capture more realistic scenarios. For instance, it is not fair to assume that everyone you come in contact with has an equal chance to infect you. A family member or a spouse is much more likely to infect you than a casual acquaintance. To capture these heterogeneous effects, we develop heterogeneous network models next.

\subsubsection{Heterogeneous network models}
%

Here we describe the dynamics of the SIS model with heterogeneous recovery and infection rates
over arbitrary strongly connected directed graphs~$\GG = (V,E)$. The recovery rate of node $i$ is given by $\delta_i > 0$. In our exposition, we consider an edge-dependent infection rate. In other words, the infection rate at which a node~$i$ is affected by an infected node~$j$ is given by $\beta_{ij}>0$ if $(i,j) \in E$.
For simplicity, we let $\beta_{ij} = 0$ if $(i,j) \notin E$. The dynamics of the SIS model in an arbitrary network is described by~\cite{AL-JAY:76}
\begin{align}\label{eq:network_sis_hetero}
\dot{p}_i &= - \delta_i p_i + \sum_{j=1}^N \beta_{ij} p_j (1 - p_i),
\end{align}
where $p_i \in [0,1]$ can be seen as either the fraction of the $i$-th subpopulation that is infected (in the metapopulation case),
or the probability that an individual $i$ is 
infected~\cite{AL-JAY:76,EC-PVM:12,CL-RB-PVM:12,AK-TB-BG:14,AK-TB-BG:14b,PVM-JO:14}.

In this model, the disease-free equilibrium $p_i = 0$ for all $i \in \until{N}$ 
is again a trivial equilibrium. In what follows, we derive conditions for this equilibrium
to be globally asymptotically stable. Let $p = (p_1, \dots, p_N)^T$ denote the state vector of the system, $D = \diag{\delta_1, \dots, \delta_N}$ the diagonal matrix of recovery rates, 
and $B = [\beta_{ij}]$ the matrix of infection rates. The dynamics~\eqref{eq:network_sis_hetero} can then be written as
\begin{align*}
\dot{p} = (B-D)p + h,
\end{align*}
where $h_i = -\sum_{j=1}^N \beta_{ij} p_i p_j$. 
The following result from~\cite{AL-JAY:76,AF-AI-GS-JJT:07,AK-TB-BG:14} characterizes
the convergence properties of these dynamics.

\begin{theorem}[Threshold condition for heterogeneous networks]\label{th:network_hetero_threshold} 
Given the dynamics in~\eqref{eq:network_sis_hetero}, for any $p(0) \neq 0$, the equilibrium $p^* = 0$ is
globally asymptotically stable if and only if $\lambda_\text{max}(B-D) \leq 0$. Furthermore,
for $\lambda_\text{max}(B-D) > 0$, there exists $p^{**} \in \real_{(0,1)}^N$ such that
$p^{**}$ is globally asymptotically stable. 
\end{theorem}

These stability results have been recently extended to a number of more complicated models such as the three-state SAIS model~\cite{FDS-CS:11}, the four-state G-SEIV model~\cite{CN-VMP-GJP:14-cdc}, and even the $SI^*V^*$ model with an arbitrary number of states~\cite{CN-MO-VMP-GJP:15}.

Now that we have provided a basic understanding of how the SIS process 
evolves and the connections between the stochastic processes
and their deterministic approximations as discussed in 
Remark~\ref{re:network}, we are now ready to formulate
and study some relevant control problems.

%
%
%
%
%
%
%

\section{Control of Epidemics}

In the previous section we presented several approaches for modeling the dynamics of spreading processes taking place on arbitrary contact networks. We have also analyzed these models and introduced several stability results for both the deterministic and stochastic cases. In this section, we describe several results aimed at controlling the dynamics of the spreading processes. 

Ideally, we are interested in controlling the stochastic network models to stop the spreading
of a disease as quickly as possible. However, before getting to the details, we must begin by
talking about our effective `control levers' in treating an epidemic. For simplicity, let us consider
the heterogeneous SIS dynamics~\eqref{eq:network_sis_hetero} 
\begin{align*}
\dot{p}_i = - \delta_i p_i + \sum_{j=1}^M \beta_{ij} p_j (1 - p_i),
\end{align*}
as a meta-population model with $M$ subpopulations. That is, 
each node~$i$ is some subpopulation (such as a town) of $n_i$ individuals in a 
larger population (such as a country) 
of $N$ individuals (see ``Meta-Population Models'' for further details). 
The parameters we have to play with are then the recovery rates $\delta_i$ for each
subpopulation and the infection rates $\beta_{ij}$ that describe the interactions
between various subpopulations.

In order to mitigate the effects of an epidemic, in general we would like to
increase the recovery rates~$\delta_i$ and decrease the infection rates~$\beta_{ij}$. 
Increasing the recovery rate of a given subpopulation can be done by providing better
treatment to sick individuals. For instance, by allocating more resources to this
particular subpopulation they can afford more doctors or better methods of treatment
for fighting a particular disease. Decreasing infection rates can be done in numerous
ways. Limiting traffic/travel between subpopulations can help decrease the infection rate.
Completely quarantining a subpopulation~$i$ is equivalent to setting~$\beta_{ji} = 0$
for all $j$ since~$i$ can no longer affect other subpopulations. Other ways of decreasing
infection rates include milder methods of prevention, such as distributing masks to
a population to minimize chance of infection; or even simply raising awareness about
a disease to make people less likely to contract the disease. 

Clearly if we had infinite resources, treatment power, and simply quarantined everyone
the disease would likely die quickly; however, this is not a feasible solution. Thus,
given a fixed budget of some sort, it is imperative to identify which parameters specifically
are the most important in order to mitigate the effects of the disease as much as possible.
We formulate these problems and discuss the current state of the art next.

\subsection{Spectral control and optimization}

Here we are interested in various optimal resource allocation problems. More specifically,
given a fixed budget, the idea is to optimally invest resources to best hinder the
spreading of a disease. Leveraging the results of Theorems~\ref{th:SIS-threshold-homogeneous}-\ref{th:network_homo_threshold}, a natural option to mitigate the effects of a possible epidemic is to make $\lambda_\text{max}(B-D)$ as small as possible. 

We first discuss the homogeneous SIS dynamics~\eqref{eq:network_sis} where $\delta$ and $\beta$ are fixed parameters. Hence, we are interested in making $\lambda_\text{max}(A)$ as small as possible. This can be achieved by modifying the network structure. The work~\cite{VMP-AJ:09} studies the effects that the network structure has on this maximum eigenvalue. In this work, the authors study how to decrease $\lambda_\text{max}(A)$ in one of two ways.
The first is to remove nodes from $A$. This might physically be done by
either quarantining or immunizing certain individuals, 
making them unable to contract the disease and, more importantly,
to spread it. Another way to reduce~$\lambda_\text{max}(A)$ is to remove links rather than completely 
removing nodes. This might physically be done by limiting traffic between certain cities
or limiting interactions between certain individuals. The caveat is that we are interested in doing this while removing the least amount of nodes or edges since these actions are likely quite costly in the real world.
The node and link removal problems of interest are then described as follows. 

\begin{problem}[Optimal node removal]\label{prob:node}
Given an original graph~$A$ and a fixed budget $C>0$, minimize~$\lambda_\text{max}(A)$ 
by removing at most~$C$ nodes from~$A$.
\end{problem}

\begin{problem}[Optimal link removal]\label{prob:link}
Given an original graph~$A$ and a fixed budget $C>0$, minimize~$\lambda_\text{max}(A)$ 
by removing at most~$C$ links from~$A$.
\end{problem}

Unfortunately, the node and link
removal problems described above are NP-complete and NP-hard, respectively~\cite{PVM-DS-FK-CL-RVDB-DL-HW:11}.
As a result, several papers instead solve convex relations or propose heuristics 
to approximately solve these problems. A simple example is one in which the nodes with
the highest degrees (largest numbers of neighbors) are removed one by one until the budget is
exhausted. Other heuristics are based on various network metrics, such as betweenness centrality~\cite{PH-BJK:02}, PageRank~\cite{JCM-JMH:07}, 
or susceptible size~\cite{CMS-TM-SH-HJH:11}, to decide which nodes should be removed first.
Similarly, there are works that are concerned with link removal rather than node 
removal~\cite{DHZ-SRG:08,PVM-DS-FK-CL-RVDB-DL-HW:11,SS-AA-BAP-AKSV:15}. In~\cite{ANB-IS:11},
the authors solve a convex relaxation of the problem and effectively project its optimal
solution onto the original problem. 

Unfortunately, the authors in~\cite{MZ-VMP:14} study the worst-case scenarios of these suboptimal strategies to show that network-based heuristics can perform arbitrarily poorly. Thus, it is
hard to evaluate a priori how well a suboptimal solution to Problems~\ref{prob:node} and~\ref{prob:link} will perform. Furthermore, completely removing nodes or even links might
not be feasible solutions anyway as this would require fully quarantining certain
subpopulations or completely shutting down certain roads or methods of travel between
various subpopulations.

Instead, let us now consider the heterogeneous network model in~\eqref{eq:network_sis_hetero} and tune the values of the parameters $\delta_i$ and $\beta_{ij}$ rather than completely changing the network structure. The authors in~\cite{EG-JO-PVM:11} formulate this as an optimization problem to minimize the steady-state infection values over heterogeneous recovery rates. A gradient descent algorithm is then proposed to find feasible local minima solutions. Another alternative is to utilize the result of Theorem~\ref{th:network_hetero_threshold}. In this direction, several papers consider the minimization of $\lambda_\text{max}(B-D)$ under various constraints. The effect of minimizing this eigenvalue is to maximize the exponential decay rate of the system towards the disease-free equilibrium.

While tuning the spreading and recovery rates, one can consider a discrete optimization setup in which one can only tune these rates within a discrete set of feasible values. This problem has been shown to be NP-complete in~\cite{BAP-LA-TI-HT-CF:13}. Alternatively, one can consider a relaxation in which these rates can take values in a feasible continuous interval. In this case, the works in~\cite{YW-SR-AS:07,YW-SR-AS:08} propose efficient methods for allocating resources to minimize the dominant eigenvalue of relevant matrices. In~\cite{VMP-MZ-CE-AJ-GJP:13b} and~\cite{VMP-MZ:13}, the problem of minimizing~$\lambda_\text{max}(B-D)$ is cast into a semidefinite program 
framework for undirected networks.
In~\cite{XZ-LZ-JW-CWT:13,VMP-MZ-CE-AJ-GJP:13}, this problem is solved for directed graphs using geometric programming where the solution can be obtained using standard off-the-shelf convex optimization software. Furthermore, geometric programs allow for the simultaneous optimization over both the infection rates and recovery rates. See ``Geometric Programming'' for further details.

In what follows, we present a simplified version of the optimization problem considered in~\cite{VMP-MZ-CE-AJ-GJP:13}
and show how it can be efficiently reformulated as a geometric program. Consider
the deterministic heterogeneous SIS model~\eqref{eq:network_sis_hetero} with natural recovery
rates $\delta_i = \underline{\delta}_i > 0$ and infection rates $\beta_i = \overline{\beta}_i > 0$
for all~$i \in \until{N}$, where $\beta_{ij} = \beta_i$ for $j \in \Nin_i$ and $\beta_{ij} = 0$ otherwise. In other words, the rate at which a node~$i$ is infected
is a node-dependent parameter rather than an edge-dependent one. We then assume we are able to pay
some cost to increase $\delta_i$ up to some maximum $\overline{\delta}_i > \underline{\delta}_i$. Alternatively, we can also pay a cost to decrease $\beta_i$ down to some minimum $\underline{\beta}_i < \overline{\beta}_i$. The control parameters are then given
by $\delta_i$ and $\beta_i$, where $\underline{\delta}_i \leq \delta_i \leq \overline{\delta}_i, \quad \underline{\beta}_i \leq \beta_i \leq \overline{\beta}_i $.

Assume we have access to cost functions $f_i(\delta_i)$ and $g_i(\beta_i)$ describing the associated cost to increase $\delta_i$ and decrease $\beta_i$, respectively. In this context, given a fixed budget $C > 0$, our
goal is to minimize~$\lambda_\text{max}(B-D)$ while satisfying the constraint that the total cost
does not exceed the given budget. This problem is formally stated below.

%
%
%

\begin{problem}[Budget-constrained allocation]\label{prob:budget-constrained}
Given a fixed budget $C > 0$,
{\rm
\begin{align*}
\begin{array}{cl} 
\underset{ \{\beta_i, \delta_i \}_{i=1}^N }{\text{minimize}} & \lambda_\text{max}(B-D)  \\ 
\text{such that} & \sum_{i=1}^N f_i (\beta_i) + g_i(\delta_i) \leq C \\
& \underline{\beta}_i \leq \beta_i \leq \overline{\beta}_i \\
& \underline{\delta}_i \leq \delta_i \leq \overline{\delta}_i .\end{array}
\end{align*}
}
\end{problem}

Note that solving Problem~\ref{prob:budget-constrained} is not trivial since the objective (maximum eigenvalue)
function is not convex in general. However, the following result guarantees that, under mild assumptions
on the cost functions, this problem can be solved
exactly by rewriting it as a geometric program which can be solved using standard off-the-shelf
convex optimization software. See~\cite{VMP-MZ-CE-AJ-GJP:13} for further details on this equivalence.

%

\begin{theorem}[Solution to budget-constrained allocation problem]
Problem~\ref{prob:budget-constrained} can be solved by solving the following auxiliary geometric program 
{\rm
\begin{align}\label{eq:budget-solution}
\begin{array}{cl}
\underset{ \lambda,\left\{ \beta_{i},\widetilde{\delta}_{i}, u_i\right\}_{i=1}^N }{\text{minimize}} & \lambda \\
\text{such that} & \sum_{j=1}^{N}a_{ij} \beta_{i}u_{j} + \widetilde{\delta}_{i}u_{i} \leq\lambda u_{i},\\ 
& \sum_{j=1}^{N} f_j(\beta_j) + \widetilde{g}_j(\widetilde{\delta}_j) \leq C,\\
 & \phi - \overline{\delta}_i \leq \tilde{\delta}_i \leq \phi - \underline{\delta}_i , \\
 & \underline{\beta}_i \leq \beta_i \leq \overline{\beta}_i, 
\end{array}
\end{align}
}
for all $i \in \until{N}$, with $\phi > \max_{j} \overline{\delta}_j$ and
$\widetilde{g}_j(\widetilde{\delta}_j) = g_j(\phi - \delta_j)$, 
where $\beta_i^*$ and $\delta_i^* = \phi - \widetilde{\delta}_i^*$ solve Problem~\ref{prob:budget-constrained} with rate $\lambda_\text{max}(Q) \leq \lambda^* - \phi$. 
\end{theorem}

We have only considered node-dependent infection rates $\beta_i$ here. This is extended to the case where
these rates can be controlled over edges $\beta_{ij}$ in~\cite{VMP-MZ-DS:14}. 
Aside from the discussed SIS model, other recent works have also applied these ideas to more general models.
The authors in~\cite{VMP-FDS-CS:13} formulate the SDP problem for a three-state SAIS model 
proposed in~\cite{FDS-CS:11} in which alertness to
a possible epidemic is also modeled. A general four-state SEIV model is considered in~\cite{CN-VMP-GJP:15-TCNS}
for which the authors develop equivalent geometric programs to optimize the dominant eigenvalue over various
parameters of the model simultaneously. 

These types of optimal allocation strategies have been recently compared to fair strategies in~\cite{AV-SR:12}, where resources must be allocated evenly across all nodes, to show their effectiveness in targeting resources rather than evenly spreading them. However, there are still
some drawbacks of these spectral control approaches that need to be addressed before
we can fully take advantage of their solutions in weakening the impact of diseases in the future.

The first main drawback of these approaches is that they do not take
into account the current state of the system. This means that even nodes that are not at
immediate risk of being infected might be allocated resources to raise their recovery rates
or decrease their infection rates. Second, solving these problems exactly requires a lot of knowledge. In addition to knowing the natural recovery rates and infection rates, we have assumed that we have exact knowledge of the entire graph which is a bit of a stretch. Third, these are centralized solutions which may take a long time to compute. Although we have been able to solve some variants of this problem as discussed above efficiently (in polynomial time), this still may not be fast enough if these networks are very large. Lastly, we have also assumed that once the optimal solution is found, we are able to instantaneously set the recovery rates and infection rates to the desired values.

We will discuss the current efforts on how each of these problems are currently being addressed and what still needs to be done in the sections to come. We begin in the next section by relaxing this first issue by looking at optimal control problems with feedback, rather than one-time optimal resource allocation solutions.

\subsection{Optimal control}

Here we discuss various optimal control problems formulated for mitigating epidemics under the SIS and SIR dynamics. However, since very little work has been done for the network models thus far, we start by looking at the classical models. 

\subsubsection{Classical models}

We begin by recalling and slightly modifying the SIS population model~\eqref{eq:population_deterministic} to account for a control action. 
Following~\cite{GAF-CAG:07}, we rewrite the original SIS population model with $\delta = \delta_1$
\begin{align}\label{eq:population_sis_natural}
\dot{p}^I = \beta p^I (1-p^I) - \delta_1 p^I ,
\end{align}
where $\delta_1 > 0$ is the natural recovery rate of an individual. We now assume that we are able to 
control this system by increasing the recovery rate of individuals in the
population from $\delta_1$ to $\delta_2 > \delta_1$. This can be achieved, for instance, by allocating antidotes or
providing other forms of treatment to a fraction of the population. Our control signal $u \in [0,1]$ is then the fraction of the population
that we provide treatment to. For simplicity, we assume that we are able to instantly
affect the recovery rates of any number of individuals in the population. The dynamics of the controlled SIS population model is then given by
\begin{align}\label{eq:population_sis_control}
\dot{p}^I = \beta p^I (1-p^I) - \left( (1-u) \delta_1 + u \delta_2 \right) p^I .
\end{align}
Applying the result of Theorem~\ref{th:population_threshold}, we obtain the following
corollary for a fixed $u(t) = \bar{u}$. 

\begin{corollary}[Population dynamics threshold condition]\label{co:population_threshold}
The solution of $p^I(t)$ approaches $0$ as $t \rightarrow \infty$ for 
\begin{align*}
\bar{u} \geq \frac{\beta - \delta_1}{\delta_2 - \delta_1}.
\end{align*}
\end{corollary}
Since $\bar{u} \in [0,1]$, Corollary~\ref{co:population_threshold} implies that
if $(\beta - \delta_1)/(\delta_2 - \delta_1) > 1$, the disease is too strong
and will never die out regardless of our control choices. On the other hand,
when $\delta_1 \geq \beta$ the natural recovery rate is high enough
to ensure extinction of the disease without any control action ($\bar{u} = 0$). Otherwise, 
it should be easy to see now that if we want to exert a minimal amount 
of a fixed control signal to ensure extinction of the disease, we can simply choose $\bar{u} = \frac{\beta-\delta_1}{\delta_2-\delta_1}$. However, it may be the case that we are still willing to use control such that the infection dies out faster than it would naturally. For instance,
having a population with many sick individuals could incur a drastic social cost that could
be instead offset by a smaller cost of treatment. We formulate this tradeoff as an optimal control problem next.

Let the cost of treatment be linear with the number of individuals treated,
and similarly, let the cost of infection be linear with the number of infected individuals.
We are then interested in minimizing the objective function
\begin{align}\label{eq:population_objective}
J_T = \int_0^T (cp^I(t) + du(t)) dt,
\end{align}
where $c > 0$ is associated with the cost of infection, $d > 0$ is the associated
with the cost of treatment, and $T > 0$ is the time horizon. 
Utilizing Pontryagin's maximum principle, it can be shown~\cite{RM-KHW:74,GAF-CAG:07,SE-MK-SS-SV:15} that the optimal solution is given by
\begin{align*}
u^*(t) \in \left\lbrace  \begin{array}{cl} \{ 0 \} & \text{for } f(t) > 0, \\
\left[ 0,1 \right] & \text{for } f(t) = 0, \\
\{ 1 \} & \text{otherwise,} \end{array} \right.
\end{align*}
with
\begin{align*}
f(t) = \psi p^I (\delta_2 - \delta_1) + d,
\end{align*}
where $\psi$ is the costate variable with dynamics
\begin{align*}
\dot{\psi} = c - \psi ( \beta(1 - 2p^I) - ( (1-u) \delta_1 + u \delta_2 ) ).
\end{align*}

It can now be shown~\cite{GAF-CAG:07} that for $\beta/(\delta_2-\delta_1) < c/d$, 
the optimal solution is to initially treat the entire population until some time
$t'$ at which nobody should be treated. For $\beta/(\delta_2-\delta_1) > c/d$, 
the optimal solution is $u(t) = 0$ for all $t \in [0,T]$.
This bang-bang solution with at most one switch is very common in similar problems.
Other works with this same kind of solution have been studied in many different variations
of this problem that consider efficiency of control~\cite{KHW:75} or control
over both $\delta$ and $\beta$ simultaneously~\cite{EH-TD:11}. Other models
have also been considered such as the SIR model~\cite{GZ-YHK-IHJ:08}
with different incidence rates~\cite{HB:01,TKK-AB:11} or a four-state 
SIRD model~\cite{MHRK-SSV-SS:11}.

Although the bang-bang solution is common, it is possible to obtain other types of solutions for different formulations of the optimal control problem. For example, it is shown in~\cite{MHRK-SS-EA:10,MHRK-SS-EA:11} that for alternative problem formulations, the optimal solution may not be a bang-bang controller for certain classes of cost functions. In~\cite{DI-GL:08}, an SIR model with quadratic control costs over both $\delta$ and $\beta$ is considered. In this case, the optimal solution is again not a bang-bang controller. A four-state SIRC model for which the optimal solution is again not a bang-bang controller is considered in~\cite{DI-NS:13}.

In the problem above, we have assumed that we are able to 
change the control signal instantaneously. Other works
consider the case in which the rate of the control signal (its time derivative) can be controlled 
instead~\cite{AA:73,AA:74,EV-FD-ME:05,PC-SC-KC:14}. We omit the technical details
of these works as the methods are very similar to the example presented above.
Interestingly, the results from these works often admit bang-bang controllers with at most one switch as optimal solutions as well. As a final note we acknowledge that in different contexts one may be interested in the problem of maximizing the impact of a spreading process (for instance a viral marketing campaign)~\cite{TKK-AB:11,AK-PD:12} rather than minimizing it.

\subsubsection{Network models}

As mentioned before, the population models are quite crude in general as they lump an
entire population's state into just a few numbers. Thus, we now turn to optimal control
of networked models but note that, so far, very little work has been done in this realm.
Three relevant papers that consider this problem in the 
context of networks are~\cite{MB-TA-TB:09,AK-TB:14,SE-MK-SS-SV:15}.
We start our exposition by proposing an optimal control problem for SIS dynamics that has yet to be solved.

Recall the SIS network dynamics with heterogeneous recovery and infection rates~\eqref{eq:network_sis_hetero},
\begin{align}\label{eq:network_sis_hetero_repeat}
\dot{p}_i = -\delta_i + \sum_{j=1}^N \beta_{ij} p_j (1 - p_i). 
\end{align}
From Theorem~\ref{th:network_hetero_threshold} we know that a necessary and sufficient condition
for extinction is $\lambda_\text{max}(B-D) \leq 0$. In the previous section, we utilized
this result as a constraint to solve optimal allocation problems. Instead, we are now
interested in solving optimal control problems, where rather than solving a one-time
optimization to determine the curing rates $\delta_i$, we allow them to vary over time.

\begin{problem}[Optimal control of an SIS network]\label{prob:network_optimal_control}
Given a linear cost of infection $c_i$ and control $d_i$ for all $i \in \until{N}$,
minimize
\begin{align}\label{eq:network_optimal_control_SIS}
J_T = \int_0^T \left( \sum_{i=1}^N c_i p_i(t) + d_i \delta_i(t) \right) dt
\end{align}
subject to the dynamics~\eqref{eq:network_sis_hetero_repeat} and $\delta_i(t) \in [\underline{\delta}, \overline{\delta}]$ for some $0 < \underline{\delta} < \overline{\delta}$
for all $t \in [0,T]$.
\end{problem}

This, along with most of its variations, is currently an open problem. Variations
include problems similar to the optimal control problems for deterministic population
models discussed earlier: control over infection rates, non-instantaneous control,
different objective functions, etc. The only work we are aware of that has tackled
this problem is~\cite{AK-TB:14}, where the authors study the linearization 
of~\eqref{eq:network_sis_hetero_repeat} around the disease-free equilibrium. The
authors are able to show for the linear dynamics that the optimal solution is
a bang-bang controller with at most one switch, similar to many results obtained for
the population models.

Although Problem~\ref{prob:network_optimal_control} is still an open problem for the SIS dynamics, a closely related problem has been successfully solved in the context of containing computer viruses~\cite{MB-TA-TB:09,SE-MK-SS-SV:15}. Here we present
a simpler version of the problem originally posed in~\cite{SE-MK-SS-SV:15}. Consider the dynamics
\begin{align}\label{eq:network_sir_hetero}
\dot{p}_i^S &= - p_i^S \sum_{j =1}^N \beta_{ij} p_j^I - p_i^S p_i^R u_i , \notag \\
\dot{p}_i^I &= p_i^S \sum_{j =1}^N \beta_{ij} p_j^I - p_i^I \pi_i R_i u_i , \\
\dot{p}_i^R &= p_i^S p_i^R u_i + p_i^I \pi_i R_i u_i , \notag
\end{align}
where as before $p_i^S$ and $p_i^I$ are the fraction of a subpopulation that are
susceptible and infected, respectively. Then, $p_i^R = 1 - p_i^S - p_i^I$ is the
fraction of individuals that are `removed'. This refers to individuals that are
immune from the infection, whether this means they were vaccinated or recovered
from the disease and are no longer susceptible to it. Additionally, $u_i$
is the control that dictates the rate at which susceptible and 
infected individuals become removed. 

\begin{problem}[Optimal control of an SIR network for malware epidemics]\label{prob:network_optimal_control_SIR}
Given a linear cost of infection $c_i$, control $h_i^1$ and $h_i^2$, and benefit of
recovery $\ell_i$ for all $i \in \until{N}$, we would like to minimize
\begin{align}\label{eq:network_optimal_control_SIR}
J_T = \int_0^T \left(  \sum_{i=1}^N - \ell_i p_i^R + c_i p_i^I + p_i^R h_i^1 u_i + p_i^R (p_i^S + p_i^I) h_i^2 u_i \right) dt ,
\end{align}
subject to the dynamics~\eqref{eq:network_sir_hetero} 
and $u_i(t) \in [0, \overline{u}_i]$ for some $\overline{u}_i > 0$ for all $t \in [0,T]$.
\end{problem}

The following result follows from Pontryagin's maximum principle~\cite{SE-MK-SS-SV:15}.

\begin{theorem}[Optimal control of an SIR network for malware epidemics]\label{th:optimal_control_sir}
There exists $\tau_i \in [0,T]$ for all $i$ such that the optimal control
is given by
{\rm
\begin{align*}
u_i^*(t) = \left\lbrace \begin{array}{ccr} \overline{u}_i & \text{for} & t < \tau_i , \\
0 & \text{for} & \tau_i \leq t \leq T . \end{array} \right.
\end{align*}
}
\end{theorem}

Again, Theorem~\ref{th:optimal_control_sir} is consistent with many other optimal control solutions
for epidemics in that the optimal solution is a bang-bang controller with at most one switch. 
Given the recentness of these results, there are still lots of variations of this
work that need to be studied. Although the dynamics~\eqref{eq:network_sir_hetero} considered 
here is very similar to the
epidemic models we have discussed throughout the article, it is not immediately applicable
due to the term~$R_i u_i$. In the context of patching, $R$ is a state of nodes who have 
a `patch' and are thus immune, and so they can spread
this patch to healthy and infected nodes. However, this concept does not translate directly
to general epidemics; a sick person cannot get better by interacting with healthy people.
We expect these types of problems to be solved for epidemics in the very near future. 

%
%

In many of the problems we have discussed above, we have assumed 
that we have direct control of the infection rate~$\beta_{ij}$ and
the recovery rates~$\delta_i$. However, this simplistic
scenario assumes that we can control these parameters for the
entire population instantaneously which is unfeasible in the context of disease
spreading. In an effort to address this oversimplification,
there is a rising body of current work in which more realistic `control'  actions are explored.
We discuss these next.

\subsection{Heuristic feedback policies}

Here we briefly review various models that are used to capture possible human behaviors
or other countermeasures employed to deter the spreading of a disease. Rather than
explicitly attempting to control the SIS dynamics as described above, the works
we discuss here are essentially extensions to the SIS model for which stability
conditions are derived. The models are created by assuming various actions people
might take, and then the `closed-loop' system stability is analyzed. 
More specifically, rather than separately considering a model and control strategies,
the model and control strategies are co-developed to yield a sense of `closed-loop'
control model. For lack of better
terminology, we refer to these as heuristic feedback policies.

More specifically, many works consider various feedback strategies that determine when nodes or links should temporarily be removed~\cite{FHC:06,BB-AO-DL:08,JOK-SRW:93,DHZ-SRG:08,IBS-LBS:10,HZ-JZ-PL-MS-BW:11,GT-JYLB-JSB:13,IT-LK:12}. Closed-loop models
are then constructed for the various strategies whose stability properties can then be analyzed. These
strategies are generally based on some sort of perceived risk that individuals have of becoming infected, causing them to either remove links to infected neighbors or completely remove themselves from the network (for example, by staying home from work or becoming vaccinated). We begin by
looking at the simpler classical models, then later show how these can be extended to network
models.

\subsubsection{Classical models}

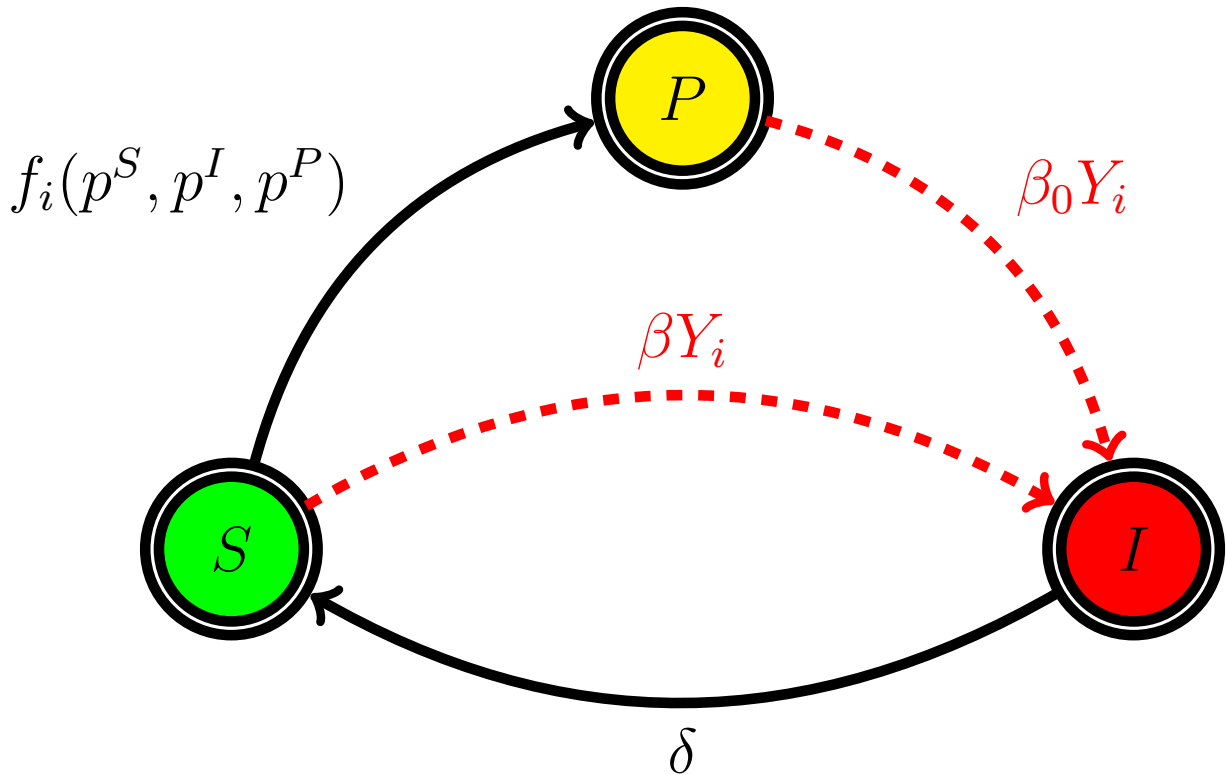
\begin{figure}[htb]
  \centering
  \scalebox{2}{\begin{tikzpicture}[shorten >= 1pt,node distance=2.25cm,auto]
  \pgfsetlinewidth{2pt}
    \node[state, fill=green, accepting] (S)              {$S$};
    \node (G) [right=of S] {};
    \node[state, fill=red, accepting] (I) [right=of G] {$I$};
    \node[state, fill=yellow, accepting] (P) [above=of G] {$P$};
    
    \path[->] (I) edge [bend left] node {$\delta$} (S);
    \path[->] (S) edge [bend left] node {$f_i(p^S,p^I,p^P)$} (P);
    
    \path[->,red,dashed] (S) edge [bend left] node{$\beta Y_i$} (I);
    \path[->,red,dashed] (P) edge [bend left] node{$\beta_0 Y_i$} (I);
    \end{tikzpicture}}
    \caption{Three-state compartmental Susceptible-Protected-Infected-Susceptible (SPIS) model. An individual~$i$
    in the infected state~$I$ transitions to the healthy or susceptible state~$S$ with
    a natural recovery rate~$\delta$. An individual in the susceptible state transitions to the
    protected state~$P$ at a rate $f_i(p^S,p^I,p^P)$ that depends on the entire network state
    and to the infected state at a rate~$\beta Y_i$ proportional to the number of infected neighbors~$Y_i$
    that node~$i$ has. An individual in the protected state transitions to the infected state at 
    a rate~$\beta_0 Y_i$ where $\beta_0 < \beta$ captures the fact that this individual is in a less
    susceptible state than normal, for instance due to behavioral changes or vaccination.}\label{fig:SAIS_compartmental}
\end{figure}

As mentioned above, these so called heuristic feedback policy solutions are all essentially
different epidemic models for which stability results are obtained. As an 
illustrative example, we consider the work~\cite{GT-JYLB-JSB:13}, where in addition to the
susceptible state~$S$ and infected state~$I$, an additional protected state~$P$
is introduced. 
The protected state refers to individuals who have decided to immunize
themselves in one way or another, and are thus not immediately susceptible to contracting the disease.
The model is described as follows. Letting $Y_i$ be the number of infected in-neighbors a
susceptible node~$i$ has in an appropriate graph, node~$i$ transitions from the susceptible
state~$S$ to the infected state~$I$
with rate~$\beta Y_i$. However, a node in the protected state~$P$ transitions to the
infected state with rate~$\beta_0 Y_i$ where $\beta_0 < \beta$ captures the 
decreased risk of infection due to being `protective'
or `alert'. A type of control is then to decide how susceptible
individuals transition to the protected state. Finally, as in the normal SIS model, individuals
that are infected naturally recover to the susceptible state with a natural recovery rate~$\delta$.
Figure~\ref{fig:SAIS_compartmental} shows the interactions of this three-state SPIS model.  
The authors then consider the extension of the SIS population dynamics~\eqref{eq:population_deterministic_whole_SIS} (by assuming a complete network topology meaning all individuals are equally likely to affect one another) to include the protected state given by
\begin{align}\label{eq:spis}
\dot{p}^S &= -\beta p^I p^S + \delta p^I - p^S f(p^S,p^I,p^P) + p^P g(p^S,p^I,p^P) , \notag \\
\dot{p}^I &= \beta p^I p^S - \delta p^I , \\
\dot{p}^P &= p^S f(p^S,p^I,p^P) - p^P g(p^S,p^I,p^P), \notag
\end{align}
where $f(\cdot)$ and $g(\cdot)$ are functions that determine how susceptible individuals are
protecting themselves. Recall that $p^S$ corresponds to the fraction of individuals in a population
that are in the susceptible state with $p^I$ and $p^P$ defined similarly for the infected and protected states, respectively. As in the case
of the deterministic SIS population dynamics~\eqref{eq:population_deterministic_whole_SIS}, one
of these equations is redundant and can be removed by using $p^S + p^I + p^P = 1$ because
we assume a constant population size. We refer to this as the three-state Susceptible-Protected-Infected-Susceptible (SPIS) model. 

The authors then explore different strategies for designing
$f$ and $g$, and analyze the stability of the system for these choices. As mentioned above,
we refer to this as a `heuristic feedback policy' because a specific control structure is
already defined and built into the model, 
rather than the objective of the work to be designing the controller itself.
More specifically, in the example above, if we are free to choose the functions $f(\cdot)$ and
$g(\cdot)$ arbitrarily, it is clear that the best thing to do is simply set $g(\cdot) = 0$
and have $f(\cdot)$ be very large. This means everybody simply protects themselves 
very quickly, for which it is easy to imagine that the disease will die out quickly as well.
Instead, it is useful to explicitly model a cost for infection and/or control as we did in the previous section. Next, we show how these models can be extended to network models.

\subsubsection{Network models}

Similar to the classical models, here we discuss models that are used to capture possible human behaviors or other countermeasures employed to deter the spreading of a disease on networks~\cite{SF-EG-CW-VAAJ:09, SF-MS-VAAJ:10}. As before, many works consider various feedback strategies that determine when nodes or links should temporarily be removed~\cite{ZR-MT-ZL:12,QW-XF-ZJ-MS:15,XLP-XJX-XF-TZ:13,CE-VMP-GJP:13,BAP-DC-NCV-MF-CF:12}. Closed-loop models are then constructed for the various strategies whose stability properties can then be analyzed. These strategies are generally based on some sort of perceived risk that individuals have of becoming infected,
causing them to either remove links to infected neighbors or completely remove themselves from the
network (for example, by staying home or becoming vaccinated).

As an example, consider again the three-state SPIS model~\eqref{eq:spis} presented in~\cite{GT-JYLB-JSB:13}.
However, we are now interested in the network version of the population dynamics there. 
To do this, we create a model very similar to the three-state SAIS model presented 
in~\cite{FDS-CS:11} where the authors introduce an alert state~$A$ which
is similar to the protected state~$P$ we consider. This state captures the possibility of human
behaviors and actions lowering the chance of contracting a disease. For simplicity, we consider
homogeneous parameters so the recovery and infection rates are set the same at all nodes. 
The deterministic version of this model is then given by
\begin{align*}
\dot{p}^S_i &=  -\beta p_i^S \sum_{j=1}^N a_{ij} p_j^I + \delta p_i^I - p^S_i f_i(p^S, p^I, p^P) ,\\
\dot{p}^I_i &= \beta p^S_i \sum_{j=1}^N a_{ij} p^I_j + \beta_0 p_i^P \sum_{j=1}^N a_{ij} p^I_j - \delta p_i^I, \\
\dot{p}^P_i &= p^S_i f_i(p^S, p^I, p^P) - \beta_0 p_i^P \sum_{j=1}^N a_{ij} p^I_j ,
\end{align*}
where $f_i(p^S, p^I, p^P)$ is a function that determines how susceptible individuals
are protecting themselves. Conditions can then be derived for the parameters
and the function~$f$ such that the disease-free equilibrium is globally asymptotically
stable~\cite{FDS-CS:11}. The authors in~\cite{FDS-CMS:12} then treat the design of this
function~$f$ as an optimal information dissemination problem. However, as in the population
dynamics case, these are very structured methods of `control' that ultimately get
built into the models. 

A very large shortcoming of these types of solutions is that they are too specific. A
very specific model with a specific control structure is proposed and studied. Unfortunately,
it is often unclear what type of spreading process each model is good for describing, if any.
In~\cite{HWH:94}, Hethcote does a great job highlighting this fact that there are far too many
slight variations of existing models. We close this section with a small anecdote 
from the epilogue of~\cite{HWH:94}, after effectively proposing $10^5$ different models, to make us really think which models are actually interesting:

\begin{displayquote}
``In the book [sic], \emph{A Thousand and One Nights}, Scheherezade had to entertain King Shahriyar with a new story each evening in order to avoid being killed. If they were mathematical biologists and she had only to present one new epidemiological model each night to entertain him, then she could have survived each night for at least 270 years. Of course, the King would probably have become disenchanted by the ``new'' models if they were only very slight variations on previous models and would have killed Scheherezade. Similarly, referees (the Kings) might become disenchanted if the papers which they receive contain models which are only slight variations on previous models. Thus I suggest that we as modelers and mathematicians should be cautious and not assume that every mathematical analysis of a slightly different model is interesting.''
\end{displayquote}

We review this issue and other technical challenges in the following section.

\section{Future Outlook}

In the previous section we have provided a high level overview of the current state of the art
involving the control of epidemics. However, there are still substantial shortcomings of
the results presented that need to be taken account for to take full advantage of their proposed
solutions. Here we highlight several of the main research challenges, how they are currently being
addressed, and what still needs to be done.

\newcounter{list}
\stepcounter{list}
\arabic{list}) All control methods we have discussed so far have been for deterministic models. 

\noindent Ideally, we are interested in ultimately controlling the original stochastic epidemic models
from which the deterministic models are derived. Results like Theorems~\ref{th:SIS-threshold-homogeneous} and~\ref{th:SIS-threshold-homogeneous-upper} help draw connections between the two,
but these have only been done for simple cases so far. Furthermore, while these results help justify
how using spectral control and optimization methods for deterministic models translate to
the original stochastic models, it is unclear how the optimal control solutions found
for deterministic models relate to the stochastic ones. 

\noindent We are only aware of a few works that attempt to control epidemic processes on networks for
the original stochastic models. In~\cite{CB-JC-AG-AS:10}, the authors consider the SIS model
with a simple heuristic control law where the curing rates~$\delta_i$ for each node are proportional to the number of neighbors. The authors are then able to show that on
any graph with bounded degree, this policy can achieve sublinear expected time to extinction
with a budget proportional to the number of nodes in the network.  A drawback of the above method is that it does not take into account the current state
of the system. This means that nodes with no infected neighbors may be assigned high
curing rates. Instead, the authors in~\cite{FC-PH-AT:09} use a heuristic PageRank algorithm 
to allocate curing resources given a fixed budget, based on the initial condition of infected
nodes. The authors are then able to provide probabilistic upper bounds on the 
expected extinction time. More recently, the work~\cite{KD-AO-JNT:14} proposes
an algorithm for which the expected time to extinction is sublinear using only a 
sublinear budget (in the number of nodes) for graphs satisfying certain technical conditions.
In~\cite{RS-NP-NF:12}, the authors consider a similar problem for which various algorithms
are developed using a Markov decision process framework.

\stepcounter{list}
\arabic{list}) All control methods we have discussed so far have admitted centralized solutions.

\noindent This is a big problem since human contact networks can be massive in practice and it may not be computationally possible to solve these problems in a centralized setting. In this direction, distributed allocation and control strategies are an interesting alternative. Again, there are only a few recent works that have looked at this problem~\cite{ER-SM:14,SE-MK-SS-SV:15,CE-AJ-VMP-GJP:15,NJW-CN-VMP-GJP:15-cdc}. As more work in optimization and control of epidemic processes is being done, we desire and expect distributed versions of these algorithms to follow.

\stepcounter{list}
\arabic{list}) All problems and solutions discussed so far assume no uncertainties.

\noindent This is a very big issue in the context of epidemics. Throughout all the modeling,
analysis, and control solutions we have presented so far, we have assumed perfect knowledge
of everything including recovery rates, state information, and network structures.
These are clear oversimplifications since in practice we would be more than lucky to have
any of these parameters simply handed to us. A review of analysis and approximation techniques considering uncertainties in the spreading parameters is provided in~\cite{LZ-GV:09}, which
is a large field of study in epidemiology.

\noindent In the context of control, far less work has been done for the case where the topology
is unknown. In~\cite{IT-LK:15}, the authors use observed infection data of a discrete time
SIS process to estimate the network topology. Optimization and control methods can then be
applied to the estimated topology; unfortunately, it is unclear how well these solutions will
perform on the actual topology. Instead, a data-driven approach to optimally allocate resources has been recently proposed in~\cite{SH-VMP-CN-GJP:15-TNSE}, where only empirical data about the spreading of a disease is available. In this work, the authors assume that the spreading and recovery rates are unknown. Alternatively, the authors assume that the responsible health agency has access to historical data describing the evolution of the disease in a network during a relatively short period of time. In this context, the authors in \cite{SH-VMP-CN-GJP:15-TNSE} propose a robust optimization framework to allocate resources based on historical data. 

\noindent Another large issue related to assuming perfect knowledge is assuming that we are able to set recovery and infection rates to whatever desired values we like. This is again a clear oversimplification and studies are needed into how various control solutions perform when
these rates cannot be set exactly.

\noindent The assumption of being able to observe exact state data is another big issue that has received little to no attention in the context of controls. These do not apply to the spectral optimization or heuristic feedback control methods, but are certainly important for the optimal
control methods.  

\stepcounter{list}
\arabic{list}) We require much more general epidemic models.

\noindent Although there has been a lot of work on modeling in addition to the SIS and SIR dynamics we have mainly focused on throughout this article, there is still a lack of generalized models. 
More specifically, a majority of works that study spreading processes begin with a single model with a fixed number of states and interactions. Many of these models are created by first looking at empirical data of a spreading process like AIDS~\cite{AV:07} or a computer virus~\cite{RPS-AV:01b}, 
then determining what type of model and how many
states should be used to capture its behavior. Instead, few works propose much more general models
with arbitrary numbers of states or layers on which the disease can 
spread~\cite{FDS-CS-PVM:13,MJK-PR:07,CN-MO-VMP-GJP:15}. The further development, analysis, and
control of these generalized models can allow rapid prototyping of models for spreading processes
that might not even exist today, in addition to completely generalizing the myriad of specific
models available today.

\noindent All models we have discussed in this article so far only consider the spreading of a single disease or process. Extending existing models to capture multiple diseases that co-evolve in a network has recently been gaining attention~\cite{MEJN:05,WM-FB:14,XW-NCV-AP-IN-MF-CF:13, AS-DT-LK:14, XW-NV-BAP-IN-MF-CF:12, AB-BAP-RR-CF:12, BAP-AB-RR-CF:12}. In general, these diseases are assumed to be mutually exclusive, meaning an individual can only be infected with one type of infection at a time. While it is discussed here in the context of disease and epidemics, these models are more aptly used in studying belief propagation or product adoption. For instance the mutual exclusion of infections is very relevant  in competition in politics, such as Democrats
vs Republicans, or competition in a marketplace, such as iPhone vs Droid vs Galaxy.
Here we briefly present a three-state two-infection $SI_1SI_2S$ model on arbitrary networks 
studied in~\cite{FDS-CS:14} and further analyzed in~\cite{NJW-CN-VMP-GJP:15-acc}. 

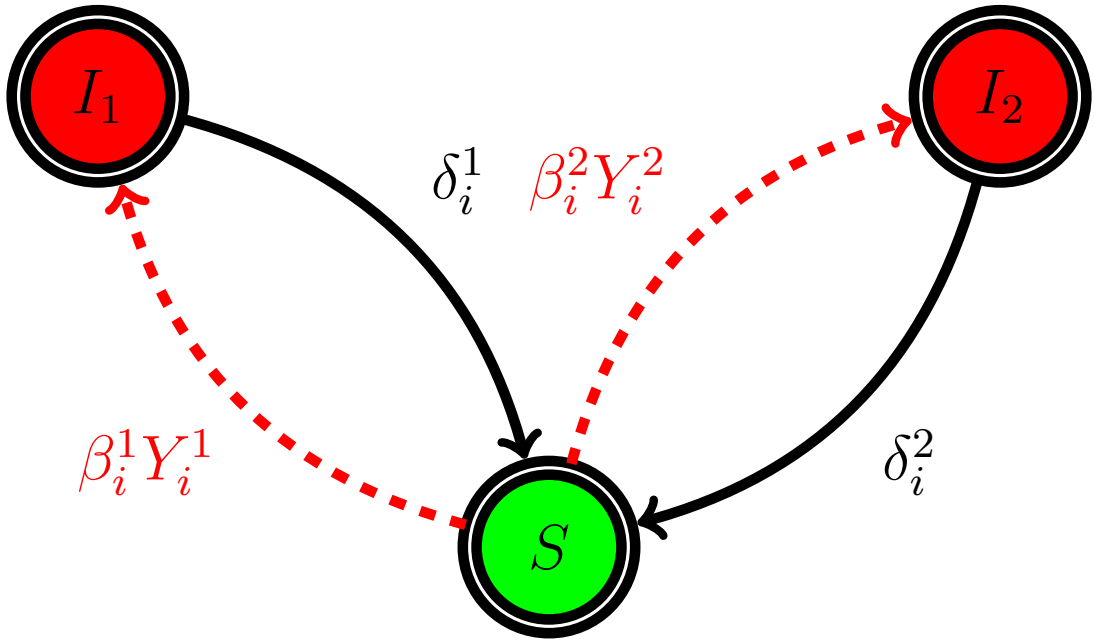
\begin{figure}[htb]
  \centering
  \scalebox{2}{\begin{tikzpicture}[shorten >= 1pt,node distance=2.25cm,auto]
  \pgfsetlinewidth{2pt}
    \node[state, fill=red, accepting] (I1)              {$I_1$};
    \node (G) [right=of I1] {};
    \node[state, fill=red, accepting] (I2) [right=of G] {$I_2$};
    \node[state, fill=green, accepting] (S) [below=of G] {$S$};
    
    \path[->] (I1) edge [bend left] node {$\delta_i^1$} (S)
    (I2) edge [bend left] node {$\delta_i^2$} (S);

    \path[->,red,dashed] (S) edge [bend left] node{$\beta^1_{i} Y^1_i$} (I1)
    (S) edge [bend left] node{$\beta^2_{i} Y^2_i$} (I2);
    \end{tikzpicture}}
    \caption{Three-state compartmental $SI_1SI_2S$ model for two diseases.
    An individual~$i$ can be in the healthy or susceptible state~$S$, or infected
    by one of two possible infections, but not both simultaneously. 
    An individual in the first infectious state~$I_1$ recovers to the susceptible state
    at a natural recovery rate~$\delta^1_i$. Similarly, an individual~$i$ in the
    second infectious state~$I_2$ recovers at a natural recovery rate $\delta^2_i$. 
    An individual in the susceptible state transitions to the infectious state$I_k$
    at rate $\beta_i^k Y_i^k$ for $k \in \{1, 2 \}$, where $Y_i^k$ is the number of
    neighbors of~$i$ that are in infectious state~$I_k$. The parameter~$\beta_i^k$
    captures the effect that neighbors of node~$i$ has on it for infection~$I_k$.
    }\label{fig:SISIS_compartmental}
\end{figure}

\noindent The model is described as follows. Let $Y_i^1$ be the number of neighbors of node~$i$
infected by the first disease~$I_1$. A node~$i$ in the susceptible state~$S$ transitions 
to the infected state~$I_1$ with rate~$\beta_i^1 Y_i^1$. Similarly, a node~$i$ in the susceptible
state transitions to the infected state~$I_2$ with rate~$\beta_i^2 Y_i^2$. Each node has its
own recovery rate for each disease given by $\delta_i^1$ and~$\delta_i^2$. For example,
a node~$i$ in the infectious state~$I_1$ recovers to the susceptible state at rate~$\delta^1_i$.
Figure~\ref{fig:SISIS_compartmental} shows the interactions of this three-state $SI_1SI_2S$ model.  
The deterministic version of this model is then given by
\begin{align*}
\dot{p}_i^S &= -p_i^S \sum_{j=1}^N a_{ij} (\beta_i^1 p_j^{I_1} + \beta_i^2 p_j^{I_2}) +  \delta_i^1 p_i^{I_1}+ \delta_i^2 p_i^{I_2}   ,\\
\dot{p}_i^{I_1} &= p_i^S \sum_{j=1}^N a_{ij} \beta_i^1 p_j^{I_1} - \delta_i^1 p_i^{I_1}, \\
\dot{p}_i^{I_2} &= p_i^S \sum_{j=1}^N a_{ij} \beta_i^2 p_j^{I_2} - \delta_i^2 p_i^{I_2} .
\end{align*}
For simplicity, we have only presented this model assuming both infections evolve over the same graph structure~$A$.
Instead, the works~\cite{FDS-CS:14,NJW-CN-VMP-GJP:15-acc} provide analysis for these dynamics over possibly different structures. A few recent works have studied the problem of controlling multiple diseases in different scenarios~\cite{GE-QZ:13,XC-VMP:14,NJW-CN-VMP-GJP:15-acc}; however, these works are still in their infancy and there are still many open problems left to be solved.

\noindent All the works about epidemics on networks we have discussed in this article have assumed a fixed graph structure. However, it is easy to surmise that
this may not be a good assumption depending on the time-scale of a spreading process. For instance,
in the context of diseases, the network of contacts in a human population is constantly changing. Hence, a time-varying network model might be more appropriate, albeit more challenging to analyze. There is still very little work analyzing these types of time-varying models. Some recent works exist~\cite{EHV-PC-RM-FB:11,SX-WL-LX-ZZ:14,MO-VMP:15} 
that have begun tackling these problems and laying the foundation for future works in this branch
of epidemic research. As with optimal control, similar problems have been studied in different contexts such as information dissemination in mobile networks~\cite{MS-KRM-BK-FB:14}, but far less has been considered in the context of epidemics thus far.

\noindent In addition to the models we have presented throughout the article, it is
worth mentioning that many works present
the same types of models from a game-theoretic 
perspective~\cite{CTB-DJDE:04,SF-EG-CW-VAAJ:09,SF-MS-VAAJ:10,AK-TB:14b,JM-AO-PVM:09}.
This is another space in which there is not yet a significant amount of study, but some seminal
works have shown its usefulness in modeling spreading processes, especially in the context of control
and optimization~\cite{YH-ST-EA-HW-PVM:14,ST-YH-EA-HW-PVM:14,AF-AJ:12,AF-AA-AJ:14}.

\section{Conclusions}

This article has reviewed and analyzed some of the most popular models studied in epidemiology. In particular, we have presented deterministic and stochastic models in the context of both population and networked dynamics. We have described many results concerning the optimization and control of epidemic dynamics, while also outlining a number of new avenues for further exploration in this field. Although the focus of this article was on disease and epidemics, it should be emphasized that the same mathematical tools
and results apply almost directly to a vast number of other spreading processes including
information propagating through a social network, malware spreading in
the World Wide Web, or viral marketing.

Despite the vast literature studying the problems discussed in this article, there are many interesting control problems left to be solved, particularly those in the context of networked dynamics. There is plenty of work left to be done to really harness the power of these results and make
a real societal impact; especially in understanding how to effectively control these
processes on complex networks. In this respect, control engineers truly have a lot to  offer in this reemerging field of research.

%
%
%
%
%
%
%
%

%
%
%
%
%
%

%
%
%
%

\clearpage



%
%

\clearpage

\clearpage

\clearpage

\clearpage

\clearpage

\clearpage

\clearpage
%
%

\setcounter{equation}{0}
\renewcommand{\theequation}{S\arabic{equation}}
\setcounter{table}{0}
\renewcommand{\thetable}{S\arabic{table}}
\setcounter{figure}{0}
\renewcommand{\thefigure}{S\arabic{figure}}

\section{Sidebar 1 \\ Graph Theory}

A \emph{graph} is a mathematical description of a given network. A graph
consists of different nodes, or vertices, and links between the nodes,
or edges, that describe the interactions between the nodes. In the context
of epidemics, a the meaning of a single node depends on the granularity
of the considered model. For example, a \emph{node} at the lowest level can
represent a single person and \emph{links} to other nodes can represent the
interactions this person has with others. On a much higher level, a single
node can represent an entire city of people, and links to other nodes can
represent the interactions this city has with others; for example, traffic flow
between cities. See ``Meta-Population Models'' for further details.

Formally, we define a \emph{directed graph}~$\GG = (\VV, \EE)$ as a pair consisting of a set
of~$N$ vertices~$\VV$ and a set of edges~$\EE \subset \VV \times \VV$.
The adjacency matrix~$A \in \realnonnegative^{N \times N}$ of~$\GG$ satisfies
$a_{ij}=1$ if and only if $(v_i,v_j) \in \EE$.  Edges are directed,
meaning that they are traversable in one direction only. The sets of
\emph{in-neighbors} and \emph{out-neighbors} of $v \in \VV$ are
respectively
\begin{align*}
  \Nin(v) &= \setdef{ v' \in \VV }{ (v',v) \in \EE} , \\
  \Nout(v) &= \setdef{ v' \in \VV }{ (v,v') \in \EE}.
\end{align*}
We say that a graph is \emph{undirected} if for all $a_{ij} = 1$, it is
also true that $a_{ji} = 1$. In this case the set of in-neighbors and
out-neighbors for each node are identical.

A \emph{directed path} $P$, or in short path, is an ordered
sequence of vertices such that any two consecutive vertices in $P$
form an edge in $\EE$. A graph~$\GG$ is \emph{strongly connected}
if for all vertices $v \in \VV$, there exists a path to all other 
vertices $v' \in \VV$.

\clearpage

\setcounter{equation}{0}
\renewcommand{\theequation}{S\arabic{equation}}
\setcounter{table}{0}
\renewcommand{\thetable}{S\arabic{table}}
\setcounter{figure}{0}
\renewcommand{\thefigure}{S\arabic{figure}}

\section{Sidebar 2 \\ Meta-Population Models}

Throughout this article we often refer to individuals and discuss the state of
all individuals in a network. However, especially in the context of diseases
spreading through populations, the number of individuals~$N$ in a given network
can be quite large. Instead of considering the entire population of interest
together, meta-population allow groups of individuals to be lumped together
into subpopulations under some assumptions.

Consider the heterogeneous network SIS dynamics~\eqref{eq:network_sis_hetero}
\begin{align}\label{eq:network_sis_sbar}
\dot{p}_i &= - \delta_i p_i + \sum_{j=1}^N \beta_{ji} p_j (1 - p_i).
\end{align}
This was essentially introduced with $p_i$ referring to the probability that
an individual~$i$ is infected (see ``Networked Mean-Field Approximations'' 
for further details). This means an $N$-dimensional system must be analyzed to
properly study how this model evolves which can be difficult if $N$ is very
large. 

Instead of studying the state of each individual in the population separately,
we can instead create $M << N$ subpopulations to approximate the dynamics of the
entire~$N$ dimensional system. 
We can then keep track of the state of each subpopulation
rather than the state of each individual in the population. 
This was originally done and analyzed 
for~$M = 2$ and turned out to be easily extendable~\cite{NTJB:86}.

Let $i \in \until{M}$ denote the $i$th subpopulation with $n_i$ individuals,
where each individual from the original population with~$N$ people is assigned
to exactly one subpopulation. In other words, we have $\sum_{j=1}^M n_j = N$. 
Note that the number of individuals in each subpopulation do not need to be the same.

We now define the dynamics of the meta-population model assuming that each subpopulation~$i$
is well-mixed and have homogeneous recovery rates~$\delta_i'$.
In other words, within each subpopulation, each individual is assumed
to have equal contact with everyone else. This is the same way the deterministic SIS
population dynamics~\eqref{eq:population_deterministic} are derived; however, we
have the extra consideration that subpopulations can affect each other as well. 
The infection rate~$\beta_{ji}'$ captures the effect that subpopulation~$j$ has
on subpopulation~$i$. Note that it is not required that~$\beta_{ji}' = \beta_{ij}'$
nor does it make sense to. Since subpopulations can have different numbers of people,
it is reasonable to think that one subpopulation~$i$ can affect another subpopulation~$j$ 
more than~$j$ can affect~$i$.
Letting $x_i$ denote the fraction of individuals in subpopulation~$i$
that are infected, we can then write the dynamics as
\begin{align}\label{eq:meta}
\dot{x}_i &= - \delta_i' x_i + \sum_{j=1}^M \beta_{ji}' x_j (1 - x_i).
\end{align}

This has now reduced the original~$N$-dimensional system into an~$M$-dimensional one.
Furthermore, it might often make more sense to consider a meta-population model instead
of an entire population one to begin with. To properly define the full network SIS
dynamics~\eqref{eq:network_sis_sbar}, we require parameters that describe
the natural recovery rates and interconnections of all individuals within the population.
Instead, it is much more reasonable to believe these parameters can be estimated 
for entire groups of people and a reasonable meta-population model can
be described with the same level of granularity. Furthermore, state information
in the meta-population model can be estimated by looking at numbers of 
infected individuals in a given subpopulation compared to the total numbers
of individuals~$n_i$ in this subpopulation.
For example, a node~$i \in \until{M}$ at the lowest level
of granularity simply recovers the full network SIS dynamics where each
node represents a single person and links to other nodes represents the
interactions this person has with others. On a much higher level, a single
node can represent an entire city of people, and links to other nodes can
represent the interactions this city has with others; for example, traffic flow
between cities. 

%
%
%

\setcounter{equation}{0}
\renewcommand{\theequation}{S\arabic{equation}}
\setcounter{table}{0}
\renewcommand{\thetable}{S\arabic{table}}
\setcounter{figure}{0}
\renewcommand{\thefigure}{S\arabic{figure}}

\clearpage

\section{Sidebar 3 \\ Geometric Programming}

Let $\mathbf{x} \in \realpositive^N$, where
$x_{1},\ldots,x_{N}>0$ denote $N$ decision variables.
In the context of geometric programs, a \emph{monomial} function $h(\mathbf{x})$ is 
a real-valued function of the form $h(\mathbf{x}) = c_0 x_{1}^{a_{1}}x_{2}^{a_{2}}\ldots x_{N}^{a_{N}}$
with $c_0>0$ and $a_{i}\in \real$ for all $i \in \until{N}$. 
A \emph{posynomial} function $q(\mathbf{x})$ is
a real-valued function that is the sum of monomials, $q({\mathbf{x}}) = \sum_{k=1}^{K}c_{k}x_{1}^{a_{1,k}}x_{2}^{a_{2,k}}\ldots x_{N}^{a_{N,k}}$,
where $c_{k}>0$ and $a_{i,k} \in \real$ for all $i \in \until{N}$ and $k \in \until{K}$.

Before stating the definition of a geometric program, the following class
of functions will be useful.
\begin{definition}
{\rm
A function $f:\real^{N}\to \real$ is \emph{convex in log-scale}
if the function 
\begin{equation}
F\left({\mathbf{x}}\right) = \log f\left(\exp {\mathbf{x}}\right),\label{eq:Convex in Log Scale}
\end{equation}
is convex in $\mathbf{x}$ (where $\exp \mathbf{x}$ indicates
component-wise exponentiation).
}
\end{definition}
\begin{remark}
Note that posynomials (hence, also monomials) are convex in log-scale
\cite{SB-LV:04}.
\end{remark}

A geometric program (GP) is an optimization problem of the form
\begin{align}\label{eq:General GP}
\begin{array}{cc}
\underset{\mathbf{x}}{\text{minimize}} & f(\mathbf{x})\\
\text{such that} & q_{i}(\mathbf{x})\leq 1,\: i=1,\dots,m, \\
 & h_{i}(\mathbf{x})=1,\: i=1,\dots,p,
\end{array}
\end{align}
where $f$ is a funcion that is convex in log-scale, $q_{i}$ are 
posynomial functions, and $h_{i}$ are monomial functions
for $i \in \until{N}$.
A comprehensive treatment of GPs is provided in~\cite{SB-SJK-LV-AH:07}.
A GP is a quasiconvex optimization problem \cite{SB-LV:04} that can
be transformed to a convex problem utilizing a
logarithmic change of variables $y_{i}=\log x_{i}$, and a logarithmic
transformation of the objective and constraint functions.
The GP in~\eqref{eq:General GP} can then be written in the transformed
coordinates by
\begin{align}\label{eq:Transformed GP}
\begin{array}{cl}
\underset{ \mathbf{y} }{\text{minimize}} & F(\mathbf{y}) \\
\text{such that} &  Q_{i}\left({\mathbf{y}}\right)\leq0,\: i=1, \dots, m, \\
& {b}_{i}^{T}{\mathbf{y}}+\log d_{i}=0,\: i=1, \dots, p,
\end{array}
\end{align}
where $Q_{i}\left({\mathbf{y}}\right) = \log q_{i}(\exp {\mathbf{y}})$ and
$F\left({\mathbf{y}}\right) = \log f\left(\exp {\mathbf{y}}\right)$. Also, given that $h_{i}\left({\mathbf{x}}\right) =  d_{i}x_{1}^{b_{1,i}}x_{2}^{b_{2,i}}\ldots x_{N}^{b_{N,i}}$,
we obtain the equality constraint above, where ${b}_{i}= \left(b_{1,i},\ldots,b_{N,i}\right)$.

Since $f\left({\mathbf{x}}\right)$ is convex in log-scale, $F\left({\mathbf{y}}\right)$
is a convex function. Furthermore, since $q_{i}$ is a posynomial (and therefore
convex in log-scale), $Q_{i}$ is also a convex function. This
shows that~\eqref{eq:Transformed GP} is a convex optimization problem in standard
form and can be efficiently solved in polynomial time~\cite{SB-LV:04}.

\clearpage

\section{Sidebar 4 \\ Networked Mean-Field Approximations}

The method of going from a stochastic model to a deterministic mean-field approximation
is certainly not one that can be overlooked. The derivations of these approximations,
their accuracy, and what they can tell us about the original stochastic models is an entire
area of research all by itself. 

We discuss briefly here how one can go from the stochastic model~\eqref{eq:stochastic_network_sis}
to the deterministic one~\eqref{eq:network_sis}. Recall the stochastic model
\begin{align*}
\begin{array}{ll} X_i : 0 \rightarrow 1 & \text{with rate } \beta \sum_{j \in \NN_i} X_j ,\\
                  X_i : 1 \rightarrow 0 & \text{with rate } \delta . \end{array} 
\end{align*}
Given the entire state $X(t)$ at some time~$t$, we can write the probability of state~$i$ at a future
time $t' = t + \timestep$ for small $\timestep$ by
\begin{align*}
P(X_i(t') = 0 | X_i(t) = 1, X(t) ) &= \delta \timestep + o(\timestep) ,\\
P(X_i(t') = 1 | X_i(t) = 1, X(t) ) &= 1 - \delta \timestep + o(\timestep) , \\
P(X_i(t') = 1 | X_i(t) = 0, X(t) ) &= \beta \sum_{j \in \NN_i} X_j(t) \timestep + o(\timestep) ,\\
P(X_i(t') = 0 | X_i(t) = 0, X(t) ) &= 1 - \beta \sum_{j \in \NN_i} X_j(t) \timestep + o(\timestep) . 
\end{align*}
As we take $\timestep$ to 0 in these forward Kolmogorov equations, we can write the 
exact dynamics of the expectation by
\begin{align*}
\frac{dE[X_i]}{dt} &= - E \left[ \delta + (1 - X_i) \beta \sum_{j \in \NN_i} X_j(t) ) \right] \\
&=- E[X_i] \delta + E \left[(1 - X_i) \beta \sum_{j \in \NN_i} X_j(t) ) \right] . \\
\end{align*}
The complication now comes from the term $E[X_i X_j]$ relating the covariance of the random variables
$X_i$ and $X_j$ with their independent probabilities. The mean-field approximation~\eqref{eq:network_sis}
(and similar ones for different variations of the stochastic model) are then obtained by assuming 
that $E[X_i X_j] = E[X_i] E[X_j]$ for all $i \neq j$. In other words, we are assuming all the random
variables are independent.

However, this is clearly not true. This means that for any fixed population
with a stochastic model, the
deterministic approximations we study are just that; approximations. 
Naturally, this begs the questions of how accurate they are in describing their
stochastic counterparts.

Although we only consider an approximation of the expected
values~$p_i$, it has actually been shown that these are upper-bounds on the actual
probabilities~\cite{PVM:11,EC-PVM:12,CL-RB-PVM:12} (this is essentially done by showing
that $E[X_i X_j] \geq 0$ for all $i \neq j$). Fortunately, this has very positive
implications on attempting to control the underlying stochastic process by using the
deterministic mean-field model: By stabilizing the deterministic
approximations, we can make claims like the ones presented in Remark~\ref{re:network}. More
specifically, if we can ensure that the disease-free equilibrium of the deterministic
model is globally asymptotically stable, then the stochastic system will reach the 
disease-free absorbing state in sublinear time (with respect to the size of the network) 
in expectation.

In~\cite{PVM-RVDP:15}, the authors begin looking at how accurate the deterministic mean-field approximations
are in describing the stochastic models, rather than just guaranteeing the upper bound. This is
still an open problem for arbitrary networks. 

Furthermore, all works above only consider the SIS dynamics. Although the recent work~\cite{NAR-BH:15} provides this
type of analysis for a three-state SIRS model, rigorous analysis for more complicated
models in general are still unsolved problems.

\clearpage

\section{Author Information}

Cameron Nowzari received the Ph.D. in Engineering Sciences from the University
  of California, San Diego in December 2010 and September 2013,
  respectively. He is currently working as a Postdoctoral Research
  Associate at the University of Pennsylvania.
  He was a finalist for the Best Student Paper Award at the 2011 American
  Control Conference and received the 2012 O. Hugo Schuck Best Paper
  Award in the Theory category. His current research interests include dynamical systems and
  control, sensor networks, distributed coordination algorithms,
  robotics, applied computational geometry, event- and self-triggered
  control, Markov processes, network science, and spreading processes on networks.
  
Victor M. Preciado received the Ph.D. in Electrical Engineering and Computer Science from the Massachusetts Institute of Technology, Cambridge, in 2008. He is currently the Raj and Neera Singh Assistant Professor of Electrical and Systems Engineering at the University of Pennsylvania. He is a member of the Networked and Social Systems Engineering (NETS) program and the Warren Center for Network and Data Sciences. His current research interests include network science, dynamic systems, control theory, complexity, and convex optimization with applications in social networks, technological infrastructure, and biological systems.

George J. Pappas received the Ph.D. in Electrical Engineering and Computer Sciences from the University of California, Berkeley, in December 1998. He is currently the Joseph Moore Professor and Chair of Electrical and Systems Engineering at the University of Pennsylvania. He also holds secondary appointments in Computer and Information Sciences, and Mechanical Engineering and Applied Mechanics. He is a member of the GRASP Lab and the PRECISE Center. He currently serves as the Deputy Dean for Research in the School of Engineering and Applied Science. His current research interests include hybrid systems and control, embedded control systems, cyberphysical systems, hierarchical and distributed control systems, networked control systems, with applications to robotics, unmanned aerial vehicles, biomolecular networks, and green buildings.

\clearpage

\tableofcontents

\end{document}